\documentclass[article]{amsart}
\usepackage{lineno,hyperref}
\modulolinenumbers[5]
\usepackage{amsfonts, amssymb, amsmath, amsthm}
\usepackage{mathrsfs} % \mathscr
\usepackage{latexsym}% utilisation des symboles LaTeX pour avoir un beau LaTeX
\usepackage{enumerate}

\usepackage{url}

\usepackage{epic}
\usepackage{graphicx}

\usepackage{mathabx} 

\input amssym.def
\usepackage{amscd}
\usepackage[mathscr]{eucal}

\usepackage{mathtools}

\newfont{\cyrr}{wncyr10}
\usepackage{tikz-cd}
\usepackage[utf8]{inputenc}

\newcommand{\thmref}[1]{Theorem~\ref{#1}}

\newcommand{\lemref}[1]{Lemma~\ref{#1}}

\newcommand{\K}{{\mathbf K}}

\renewcommand{\b}{{\mathfrak{b}}}
\renewcommand{\a}{{\mathfrak{a}}}
\renewcommand{\c}{{\mathfrak{c}}}
\newcommand{\e}{{\mathfrak{e}}}

\renewcommand{\P}{{\mathfrak{p}}}

\renewcommand{\O}{{{\mathcal{O}}}}
\newcommand{\M}{{\mathbf M}}
\newcommand{\rO}{{{\rm{O} }}}

\newcommand{\N}{{\mathfrak N}}

\newcommand{\Q}{{\mathbb Q}}

\renewcommand{\d}{{\partial}}
\newcommand{\q}{{\mathfrak q}}
\newcommand{\R}{{\mathbb R}}

\newcommand{\A}{\mathcal{A}}
\newcommand{\B}{\mathcal{B}}

\newcommand{\m}{{\mathfrak m}}
\newcommand{\f}{{\mathfrak f}}

\def\li{\operatorname{Li}}

\newtheorem{thm}{Theorem}

\newtheorem{lem}[thm]{Lemma}
\newtheorem{cor}[thm]{Corollary}

\newtheorem{rmk}{Remark}[section]

\newcommand{\sumflat}{\mathop{\sum\nolimits^{\hbox to 0pt{$\flat$\hss}}}}
\newcommand{\prodflat}{\mathop{\prod\nolimits^{\hbox to 0pt{$\flat$\hss}}}}
\newcommand{\sumsharp}{\mathop{\sum\nolimits^{\hbox to 0pt{$\sharp$\hss}}}}
\newcommand{\prodsharp}{\mathop{\prod\nolimits^{\hbox to 0pt{$\sharp$\hss}}}}

\begin{document}

\title{Representing ideal classes of 
ray class groups by product of prime ideals of small size}

\author[J.-M. Deshouillers]{Jean-Marc Deshouillers}
\email{jean-marc.deshouillers@math.u-bordeaux.fr}

\author[S. Gun]{Sanoli Gun}   
\email{sanoli@imsc.res.in}

\author[O. Ramar\'e]{Olivier Ramar\'e}
\email{olivier.ramare@univ-amu.fr}

\author[J. Sivaraman]{Jyothsnaa Sivaraman}   
\email{jyothsnaa.s@gmail.com}

\address[J.-M. Deshouillers]
{Institut de Math\'{e}matiques de Bordeaux UMR 5251
Universit\'{e} de Bordeaux 
351, cours de la Lib\'{e}ration - F 33 405 TALENCE,
France.}

\address[S. Gun]   
{Institute of Mathematical Sciences, 
	HBNI, C.I.T Campus, Taramani, 
	Chennai  600 113, 
	India.}

\address[O. Ramare]
{CNRS / Institut de Math{\'e}matiques de Marseille, 
	Aix Marseille Universit{\'e}, U.M.R. 7373, 
	Site Sud, Campus de Luminy, Case 907, 
	13288 MARSEILLE Cedex 9, France. }

\address[J. Sivaraman]   
{Chennai Mathematical Institute,
     Plot No H1, SIPCOT IT Park,
     Padur PO, Siruseri 603103,
     Tamil Nadu, 
      India.}

\subjclass[2010]{Primary: 11R44, Secundary: 11R42}
\keywords{ Brun-Titchmarsh Theorem, Selberg sieve for number fields, Product of primes in a class 
of narrow ray class group}

\begin{abstract}
 We prove that, for every modulus $\q$, every class of the
 narrow ray class group $H_\q(\K)$ of an arbitrary number field $\K$
contains a product of three unramified prime ideals $\P$ of degree
one with $\N\P\le (t(\K)\N\q)^3$, where $t(\K)$ is an
explicit function of $\K$ described in~\eqref{deftK}. To achieve this
result, we first obtain a sharp explicit Brun-Titchmarsh Theorem for ray
classes and then an equally explicit improved Brun-Titchmarsh Theorem for
large subgroups of narrow ray class groups. En route, we deduce an explicit upper bound 
for the least prime ideal in a quadratic subgroup of a narrow ray class group
and also for the size of the least ideal that is a product of degree one
primes in any given class of $H_\q(\K)$.
\end{abstract}

\maketitle

%%%%%%%%%%%%%%%%%%%%%%%%%%%%%%%%%%%%%%%%%%%%%%%%%%%%%%%
%%%%%%%%%%%%%%%%%%%%%%%%%%%%%%%%%%%%%%%%%%%%%%%%%%%%%%%
\section{Introduction and statements of the Theorems}
%%%%%%%%%%%%%%%%%%%%%%%%%%%%%%%%%%%%%%%%%%%%%%%%%%%%%%%
%%%%%%%%%%%%%%%%%%%%%%%%%%%%%%%%%%%%%%%%%%%%%%%%%%%%%%%
A classical theorem of Linnik \cite{{Linnik*44a}, {Linnik*44b}} states
that there exists an absolute constant $b >0$ such that, given a positive
integer $q$ and an invertible residue class $a \bmod q$, the smallest
prime in the aforementioned residue class is at most $q^{b}$.  Several
efforts have been made to bound $b$ from above when $q$ is large.  The
best result to date is that of Xylouris
\cite{{Xylouris*11a},{Xylouris*11}} who proved that $b$ can be taken
to be $5$ provided $q$ is large enough.

Let $\K$ be a number field, $\O_\K$ be its ring
of integers and let $\q$ be an integral ideal of $\K$. Let us
 recall the definition of the (narrow) ray class group. 
 Let $I(\q)$ be the group of fractional ideals of $\K$
which are co-prime to $\q$ and $P_{\q}$ be the subgroup of $I(\q)$
consisting of principal ideals $(\alpha)$ satisfying
$v_\mathfrak{p}(\alpha -1) \ge v_\mathfrak{p} (\q)$ for all prime
ideals $\mathfrak{p}$ dividing $\q$ and $\sigma(\alpha) >0$ for all
embeddings $\sigma$ of $\K$ in $\R$.  This is the ray class group
associated to the modulus $\q \q_{\infty}$, where $\q_{\infty}$
is the set of all real  places of $\K$. However for notational convenience,
we shall simply refer to this as the narrow ray class group
of modulus $\q$. We set
$H_{\q}(\K)=I(\q)/P_{\q}$. When $\q=\O_\K$, the group $H_{\q}(\K)$ is
the usual class group in the narrow sense.

One way to generalize Linnik's problem to number fields is the following: 
Given an integral ideal $\q$ in $\O_\K$ and a class $C$ of $H_\q(\K)$, find a 
prime ideal $\P \in C$ such that $\N(\mathfrak{p}) \leq A~\N(\q)^b$.  This problem
was first considered by Fogels \cite{Fogels*64a} which was then refined by
Weiss \cite{{Weiss*80}, {Weiss*83a}} who proved the following theorem.
%%%%%%%%%%%%%%%%%
\begin{thm}{\rm (Weiss \cite{{Weiss*80},{Weiss*83a}} 1983)}\label{t1}
  There exists an effectively computable constant $a$ such that the
  following holds. In any number field $\K$ of discriminant $d_\K$,
  every element of $H_{\q}(\K)$ contains a prime ideal $\P$ of
  $K$ such that $\N\P \le (|d_\K| \N\q)^a$.
\end{thm}
%%%%%%%%%%%%%%%%%%
Here and thereafter, $\mathfrak{N}$ denotes the absolute norm of $\K$ over $\Q$.
Note that the exponent $a$ is not explicit.  

There is another generalization of Linnik's problem which was initiated by 
Lagarias-Montgomery-Odlyzko~\cite{Lagarias-Montgomery-Odlyzko*79}.
By class field theory, finding a prime ideal in a class of the ray
class group $H_{\q}(\K)$ is the same as finding a prime ideal with a
given Artin symbol in the Galois group of the corresponding ray class
field $\K_{\q}$ over the number field $\K$.  In this set-up,
Lagarias, Montgomery and Odlyzko derived an upper bound 
for the least norm of such a prime ideal in terms of $|d_{\K_\q / \Q}|^{a}$ 
with $a$ inexplicit. For explicit $a$, see the works of 
Ahn-Kwon \cite{Ahn-Kwon*19}, Kadiri-Wong \cite{KW}, Thorner-Zaman 
\cite{Thorner-Zaman*17} and Zaman \cite{Zaman*16}.  
All these bounds are however large in terms of the
exponent of the norm of the modulus $\q$. For example, if one
takes a rational prime $p$, then
the ray class field is the $p$-th cyclotomic field whose
discriminant has absolute value $p^{p-2}$.

In this article, we would like to show that if one considers a product of three 
prime ideals in place of a single prime ideal, one can find prime ideals with 
much smaller norm when compared to the norm of the modulus. 
Furthermore, we shall make  the dependence on the ground field 
$\K$ explicit.

Our notation recalled below is classical. In brief, $n_\K$, $h_\K$, $R_\K$ and
$d_\K$ are respectively the degree, the class number, the regulator
and the discriminant of a number field $\K$ while $\alpha_\K$ is the
residue of the Dedekind zeta-function at~1.  In this set up, we have
the following theorem.
%%%%%%%%%%%%%%%
\begin{thm}\label{mainthm}
Let $\K$ be a number field, $\q$ be an integral ideal of
$\K$. Set
\begin{equation}
  \label{deftK}
  t(\K)=\max\Bigl(n_{\K}^{ 48n_{\K}^3 }|d_\K|^6(R_\K h_\K)^{n_{\K}}, ~\exp(|d_\K|^{30})\Bigr).
\end{equation}
Each element of $H_{\q}(\K)$
contains a product of three degree one primes, each of norm
at most $t(\K)^3\N\q^3$. 
\end{thm}
%%%%%%%%%%%%%%%
As the case $\K=\mathbb{Q}$ has already been treated in
 \cite{Ramare-Walker*16,Ramare-Serra-Srivastav*18} by Ramar\'e
together with Serra, Srivastav and Walker, we
may assume that $\K\neq\mathbb{Q}$ and we do so hereafter.

\begin{rmk}
Three features of this result should be underlined:  the 
dependence on the field $\K$ is completely explicit in terms of 
classical invariants of the field, and even numerically
so. However, we did not strive to get small constants. The second point is 
that the exponent in $\q$ is relatively small. The
third point is more technical; the dependence in $\q$ has the form
$\N\q^3$ and not $\N\q^3$ times a power of $\log\N\q$ and this precision 
comes from a much more refined treatment. The dependence in $\K$ is most
probably dominated by the term $\exp(|d_\K|^{30})$ that comes from  a
possible real zero abnormally close to~1. We control such a zero by
Lemma~\ref{AKzero} due to  Kadiri and Wong \cite{KW} (see also \cite{Ahn-Kwon*19}).
\end{rmk}

Theorem~\ref{mainthm} relies on four ingredients of independent
interest. We first need a Brun-Titchmarsh Theorem for elements 
of $H_\q(\K)$.
%%%%%%%%%%%%%%%%%
\begin{thm}\label{bt-tri}
 Let $\b, \q$ be integral ideals with $(\b, \q) = \O_{\K}$ and  $[\b] \in H_{\q}(\K)$. 
 Then
 \begin{equation*}
 \sum_{\substack{\P \in [\b],\\ \N\P\le X}} 1
~\le~ 
 \frac{2  X}  { |H_\q(\K)| \log (\frac{X}{u(\K) \N\q})}~,
\quad u(\K)= n_{\K}^{ 48n_{\K}^3 }|d_\K|^6(R_\K h_\K)^{n_{\K}}
\end{equation*}
provided the denominator is~$>0$.
\end{thm}  
%%%%%%%%%%%%%%%%% 
This is the number field analogue of the classical Brun-Titchmarsh
Theorem for the initial interval, see for instance \cite[Theorem
2]{Montgomery-Vaughan*73} by Montgomery and Vaughan. A
precursor of this result can be found in \cite[Theorem
4]{Hinz-Lodemann*94} by Hinz and Lodemann, though without the
dependence in $\K$ and with a slightly worse upper bound.  As these
two authors, we rely on the Selberg sieve, though with an improved
treatment of the error term, see Theorem~\ref{absolute}, and on an
estimate for the number of integral ideals recalled as
Theorem~\ref{asymfinal} below.

Our second ingredient is a Brun-Titchmarsh Theorem valid for cosets,
the analogue of \cite[Theorem 1.2]{Ramare-Serra-Srivastav*18} when
$\K=\mathbb{Q}$, and which is the topic of Section~\ref{btc}. We
notice here that in Theorem~\ref{bt-tri}, we \emph{parametrized} the
class while for Theorem~\ref{bt}, we \emph{capture} elements of a
subgroup by using multiplicative characters.

Our third ingredient is less novel. 
%%%%%%%%%%%%%%%%%%%%%%%%%%%%%%%%%
\begin{thm}\label{degreeoneprime}
Every element in $H_\q(\K)$ contains an integral ideal $\a$ such that
$\N\a \le 10^{25n_\K} n_{\K}^{7n_{\K}} |d_\K|^{4}  \N\q^3$ 
and $\a$ is product of degree one primes.
\end{thm}
\begin{thm}\label{degreeoneprimebis}
Every element in $H_\q(\K)$ contains an integral ideal $\a$ such that
\begin{equation*}
 \N\a ~\le~ F_1(\q)\N\q \, \log(3F(\q))^{n_\K^2}\,  
 \log\log(B(\K)F(\q)\N\q)^2,
 \quad
 B(\K)=n_\K^{50n_\K^3} (E(\K)\sqrt{|d_\K|})^{n_\K},
\end{equation*}
and $\a$ is product of degree one primes. Here
$F(\q)=2^{r_1}h_\K \phi(\q)/ h_{\K,\q}$,  $F_1(\q)=2^{r_1}h_\K \N\q$ and 
the constant $E(\K)$ is equal to $1000 n_{\K}^{ 12n_{\K}^2 }(R_\K/|\mu_\K|)^{1/n_\K}
\bigl[\log\bigl((2n_{\K})^{4n_\K}R_{\K}/|\mu_\K| \bigr)\bigr]^{n_{\K}}$.
\end{thm}
%%%%%%%%%%%%%%%%%%%%%%%%%%%%%%%%%%
Theorem~\ref{degreeoneprime} is enough for our purpose, but our proof gives
only the bound
$\N\a\ll_{\K,\varepsilon} \N\q~|H_\q(\K)|^{2(1+ \epsilon)}$,
while $|H_\q(\K)|$ should be enough.
By using some techniques from sieve method, Theorem~\ref{degreeoneprimebis} corrects
this defect as far as the dependence in $\N\q$ is concerned, but the
dependence in $d_\K$ becomes much worse. 

Our fourth and last ingredient shows that quadratic subgroups
of~$H_\q(\K)$ contain small degree one prime ideal.
%%%%%%%%%%%%%%%%%%%%%%%%%%%%%%%%%
\begin{thm}\label{primeinkernel}
  Let $\chi$ be a quadratic character on $H_\q(\K)$. There exists a
  prime ideal $\P$ of degree one in $\K$ such that $(\P, \q) = \O_\K$,
  $\chi(\P)=1$ and
  $\N\P \le 8 \cdot (10^{31} n_\K^7)^{n_\K}  |d_\K|^4 \N\q^2$.
\end{thm}
%%%%%%%%%%%%%%%%%%%%%%%%%%%%%%%%%%
This bound is modest as far as the exponent of $\N\q$ is
concerned, but it is completely explicit. When $\K=\mathbb{Q}$, the
question has been treated by Linnik and Vinogradov
in~\cite{Vinogradov-Linnik*66}. Their better exponent comes from the
usage of the Burgess bounds and Siegel's Theorem while we only rely on convexity (or
equivalently, on a Polya-Vinogradov inequality). The exponent we get
is $3/2+\varepsilon$ for any positive $\varepsilon$.

%%%%%%%%%%%%%%%%%%%%%%%%%%%%%%%%%%%%%%%%%%%%%
 \section{Notation and Preliminaries}
 %%%%%%%%%%%%%%%%%%%%%%%%%%%%%%%%%%%%%%%%%%%%
\label{Notation}

\subsection*{Notation}
Let $\K\neq \mathbb{Q}$ be a number field with discriminant
$|d_\K|\ge3$ (by Minkowski's bound).
Also let us set $n_\K = [\K : \Q]\ge2$ and $\q$ be an (integral) ideal of $\K$. 
The number of real embeddings of $\K$ is denoted by $r_1$
whereas the number of complex ones are denoted by $2r_2$. 
The ring of integers of $\K$ is denoted by $\O_\K$,
the narrow ray class group modulo $\q$ is denoted by $H_\q(\K)$
and the (absolute) norm is denoted by $\N$. Throughout the
article $\P$ will denote a prime ideal in $\O_\K$, $p$ will
denote a rational prime number and for any integral ideals
$\a, \b$, their lcm and gcd in $\O_{\K}$ will be denoted
by $[\a,\b]$ and $(\a,\b)$ respectively. Further
an element of $H_\q(\K)$ containing an integral ideal
$\a$ will be denoted by $[\a]$.

A sum over degree one prime ideals will be denoted by $\sumflat_{\P}$,
and $\prodflat_\P$ will denote a product over degree one prime
ideals. Similarly $\sumsharp_{\P}$ and $\prodsharp_{\P}$ denotes respectively
a sum and a product over primes $\P$ that are \emph{not} of degree
one. As a generalisation, the sign $\sumflat_{\a}$ denotes a summation
over integral ideals $\a$ whose prime factors are all of degree one.

\subsection*{Smoothings}\label{smooth}

We shall work with a generic smoothing function
$w : \mathbb{R}_{\ge 0} \to \mathbb{R}_{\ge 0}$ with the following properties:
\begin{itemize}
\item 
$w(t)=0$ when $t \ge1$ and $t\le 1/10$,
\item
$w$ does not vanish uniformly and $|w(t)|\le 1$ throughout,
\item 
$w$ is at least $n_\K + 3$ times continuously differentiable,
\item
For every $m\le n_\K+2$, we have
$w^{(m)}(\frac{1}{10})=w^{(m)}(1)=0$, where 
$w^{(m)}$ denotes the $m$-th derivative of~$w$.
\end{itemize}
We will henceforth refer to this function as  `the smoothing function'.  
Its Mellin transform $\check{w}$ is defined by
\begin{equation}
  \label{defcheckw}
  \check{w}(s)=\int_0^\infty w(t) t^{s-1}dt.
\end{equation}
We show in \lemref{smoothdecay} that this analytic function
decreases at least like $1/(1+|s|)^{n_\K + 3}$ uniformly in any vertical
strip.

For applications, we  select the special function $w_0$ described below. 
Let
$$
f_k(t) 
~=~
\begin{cases}
 ( 4t(1-t) )^k  & \text{  when  }t \in [0, 1], \\
 0 &   \text{  otherwise}
 \end{cases}
$$ 
be as defined in page 348 of \cite{Ramare*17-2} with $k=n_{\K}+4$. We set
\begin{equation}
  \label{defw0}
  w_0(t)=f_{n_{\K}+4}\left(\tfrac{10}{9}(t-\tfrac1{10})\right).
\end{equation}
%%%%%%%%%%%%%%
\begin{lem}
  \label{studyw0}
  We have
  $\displaystyle
    \|w_0\|_\infty=1,\
    10\sqrt{n_{\K}}\,\check{w}_0(1)\in[2, 15],\
    \|w_0'\|_1=2,\ 
    \|w_0^{(n_\K+3)}\|_\infty\le 4(40n_\K)^{n_\K + 3}$.
\end{lem}
%%%%%%%%%%%%%%%
%%%%%%%%%%%%%%%%%%%%%%
\begin{proof}
  Indeed, by \cite[Lemma 2.2]{Ramare*17-2}, we have
  $\check{w}_0(1)=\|w_0\|_1=\frac{9}{10}\frac{2^{2k}\cdot k!^2}{(2k+1)!}$
  with $k=n_\K+4$. Applying the classical explicit Stirling's formula
  \begin{equation}
    \label{ExplicitStirling}
    n!=(n/e)^n\sqrt{2\pi n}~e^{\frac{\theta_+(n)}{12 n}},   \qquad(\theta_+(n) \in [0,1]),
  \end{equation}
  we find that
  \begin{equation*}
  \|w_0\|_1
    =\frac{9}{10}\frac{\sqrt{\pi k}}{(2k+1)}e^{\frac{\theta_+(k)}{6k} ~-~ \frac{\theta_+(2k)}{24k}}
    =\frac{\xi(n_{\K})}{10\sqrt{n_{\K}}\,},    \qquad 2\le \xi(n_{\K})\le 15.
  \end{equation*}

  Next we check that
  \begin{equation*}
    \|f_k'\|_1=2\int_0^{1/2}f'_k(t)dt=2.
  \end{equation*}
  Also, by Leibniz Formula, we find that
  \begin{equation*}
   f_k^{(k-1)}(t)=4^k\sum_{0\le \ell\le
      k-1}\binom{k-1}{\ell}\frac{k!}{(k-\ell)!}t^{k-\ell}
    (-1)^{k-1-\ell}\frac{(k-1)!}{(k-(k-1-\ell))!}(1-t)^{k-(k-1-\ell)}
  \end{equation*}
  so that
$$
 |w_0^{(k-1)}(t)| 
 ~\le~
\left(\frac{10}{9}\right)^{k-1} 4^k(k-1)! ~\sum_{0 \le \ell \le k-1}\binom{k-1}{\ell}\binom{k}{k-\ell}
~=~ \left(\frac{10}{9}\right)^{k-1} 4^k(k-1)! \binom{2k-1}{k}
$$
by Vandermonde's identity. We bound this almost trivially:
\begin{equation*}
\left(\frac{10}{9}\right)^{k-1}4^k(k-1)!\binom{2k-1}{k}
~=~
\left(\frac{10}{9}\right)^{k-1}4^k\frac{(2k-1)!}{k!}\le 4^k(2k-1)^{k-1}
 %=2^{2n_\K+6}(2n_\K+7)^{n_\K+3}
\le 4(40n_\K)^{n_\K + 3}.
  \end{equation*}
\end{proof}
%%%%%%%%%%%%%%%%%%%%%%

%%%%%%%%%%%%%%%%%%%%%%%%%%%%%%%%%%%%%%%%%%%%%%%%%%%%%%%%%
\subsection*{Explicit number of ideals in a ray class below some bound}
%%%%%%%%%%%%%%%%%%%%%%%%%%%%%%%%%%%%%%%%%%%%%%%%%%%%%%%%%

We now state the main theorem in \cite{Gun-Ramare-Sivaraman*22b} which
is required for the proofs of \thmref{bt-tri} and \thmref{degreeoneprimebis} .
%%%%%%%%%%%%%%%
\begin{thm}\label{asymfinal}
Let $\q$ be a modulus of $\K$ and $[\a]$ be an element of $H_{\q}(\K)$.
For any real number $X \ge 1$, we get
\begin{equation*}
\sum_{\b \in [\a] \atop{ \b \subseteq \O_{\K} \atop{ \N\b \le X } } } 1 
~=~ \frac{\alpha_{\K} \phi(\q)}{h_{\K,\q}}
\frac{X}{\N(\q)} 
+
\rO^* \biggl( E(\K)
  F(\q)^{\frac{1}{n_\K}}\log(3F(\q))^{n_\K}
\left(\frac{X}{\N(\q)}\right)^{1-\frac{1}{n_\K}}
~+~ n_{\K}^{8n_\K}\frac{R_\K}{|\mu_\K|} F(\q) \biggr).
\end{equation*}
where $F(\q)=2^{r_1}\frac{h_\K \phi(\q) }{h_{\K,\q}}\ge1$,
$E(\K)=1000 n_{\K}^{ 12n_{\K}^2 }(R_\K/|\mu_\K|)^{1/n_\K}
  \bigl[\log\bigl((2n_{\K})^{4n_\K}R_{\K}/|\mu_\K| \bigr)\bigr]^{n_{\K}}$
  and the notation $\rO^*$ denotes that the implied constant is less than
 or equal to $1$.
\end{thm}
%%%%%%%%%%%%%%%%%%%%

%%%%%%%%%%%%%%%%%%%%%%%%%%%%%%%%%%%%%%%%%%%%%%%%%%%%%%%%% 
\subsection*{Lower bound for the root discriminant}
%%%%%%%%%%%%%%%%%%%%%%%%%%%%%%%%%%%%%%%%%%%%%%%%%%%%%%%%%

The root discriminant is defined by $|d_\K|^{1/n_\K}$.
%%%%%%%%%%%%%%%%%%%
\begin{lem}
  \label{rootdisc}
  We have $|d_\K|^{1/n_\K}\ge \pi/2$.
\end{lem}
%%%%%%%%%%%%%%%%%%%

%%%%%%%%
\begin{proof}
  In this proof, we denote $n_\K$ by $n$.
  By Minkowski's bound, we find that 
  $$
  \rho=|d_\K|^{1/n_\K}\ge
  \frac{\pi}{4}n^2/n!^{2/n}.
  $$ 
  The explicit Stirling Formula recalled in \eqref{ExplicitStirling} yields $n!^{2/n}\le
  \frac{n^2}{e^2}(\sqrt{2\pi n}~e^{\theta_+(n)/(12n)})^{2/n}\le n^2/2$
  when $n\ge 3$ while $2!^{2/2}=2\le 2^2/2$.

  We note that, when $n\ge 5$, the Minkowski bound is superseded by the bound 
  given in~\cite[Eq. (2)]{Liang-Zassenhaus*77} by Liang and
  Zassenhaus (the quantity $V_{r_1,r_2}$ is given on line~14, page~18
  of their paper, and $\mu_n\ge1$).
\end{proof}
%%%%%%%%

\subsection*{The Dedekind zeta-function}

For $\Re s= \sigma > 1$, the Dedekind zeta-function is defined by
$$
\zeta_{\K}(s)=\sum_{\a \subseteq \O_{\K}, \atop{\a \ne (0) }} \frac{1}{\N(\a)^s}, 
$$
where $\a$ ranges over the integral ideals of $\O_{\K}$. 
It has only a simple pole at $s=1$ of residue $\alpha_{\K}$, say. 
We know from the analytic class number formula that
\begin{equation}\label{acf}
\alpha_{\K} = \frac{2^{r_1} (2\pi)^{r_2} h_{\K} R_{\K}}{|\mu_{\K}| \sqrt{|d_\K|}},
\end{equation}
where $h_{\K}, R_{\K}, d_{\K}$ are as defined in the introduction and
$\mu_{\K}$ is the group of roots of unity in $\K$. As we will see later
(\lemref{HR1}) its order in vertical strips is sufficiently moderate so that 
multiplication by $\check{w}(s)$ makes it an $L^1$-function 
on any line $\Re s=\sigma\ge -1/2$.

%%%%%%%%%%%%%%%%%%%%%%%%%%%%%%%%%%%%%%%%%%%
\subsection*{Euler-Kronecker constant}
%%%%%%%%%%%%%%%%%%%%%%%%%%%%%%%%%%%%%%%%%%%
For a number field $\K$, the Euler-Kronecker constant
is defined as
$$
\gamma_{\K} = \lim_{s \to 1} \left( \frac{\zeta'_{\K}(s)}{\zeta_{\K}(s)}+ \frac{1}{s-1}  \right).
$$
We also know that in a neighbourhood of $s=1$
$$
\zeta_\K(s) = \frac{\alpha_{\K}}{s-1} ~+~ \alpha_{\K}\gamma_{\K} ~+~ \rO(s-1),
$$
where the constant in $\rO$ depends on $\K$. The constant $\gamma_\K$
is called the `Euler-Kronecker constant' in \cite{Ihara*06} by
Ihara. We use Proposition~3, page~431 of this paper, namely the inequality
\begin{equation}
  \label{IharaIneq}
  \gamma_{\K}\ge-\tfrac12\log|d_\K|
\end{equation}
(where $\alpha_\K$ in Ihara's paper is given by (1.2.2), $\beta_\K$ by (1.2.3) and $c_\K=1$ by
(1.3.12)). The conclusion we need is also restated as (0.7) therein.

%%%%%%%%%%%%%%%%%%%%%%%%%%%%%%%%%%%%%%%%%%%
\subsection*{The narrow ray-class group}
%%%%%%%%%%%%%%%%%%%%%%%%%%%%%%%%%%%%%%%%%%%
By narrow ray class group $H_\q(\K)$, we consider those ray class groups
where the integral ideal $\q$ is completed with all real
archimedean places. We have
\begin{equation}\label{eq:4}
h_{\K, 1}=|H_1(\K)| ~{\Big\vert}~  |H_\q(\K)| ~{\Big\vert}~  2^{r_1} \phi(\q)|H_1(\K)|,
\end{equation}
where 
\begin{equation} \label{eq:6}
\phi(\q)=\N(\q)\prod_{\mathfrak p|\q}\left(1 - \frac{1}{\N(\mathfrak p)}\right)
\end{equation}
and $H_1(\K)$ denotes the narrow ray class group
corresponding to $\O_{\K}$. A good reference for this are the notes of
Sutherland \cite{Sutherland*15}.
We also have the following theorem in this context.
\begin{lem}[Lang, \cite{Lang*94}, page 127]\label{classnumber}
Let $\q$ be a modulus of $\K$, $h_{\K, \q} =   |H_\q(\K)|$ and
$r_1, h_{\K}$ be as defined earlier. Then
$h_{\K, \q} \le 2^{r_1}  \phi(\q) ~h_{\K} $.
\end{lem}

\subsection*{Characters on the narrow ray-class group}

A manner to work with the narrow ray class group $H_\q(\K)$ is to consider
its character group. When lifted to the set of all ideals, these are
characters that vanish on ideals which are not co-prime to $\q$.

An excellent reference is the report \cite{Landau*18-2} where Landau
explains in detail and refines Hecke's theory. We are only interested
in the extensions of the true characters of $H_\q(\K)$ and these
are the ones that have finite order. These are the ones that Hecke
considers, while Landau \cite[Lemma~6.34 and (6.7)]{Narkiewicz*04} considers
an extended class of characters that may have infinite orders.

In our case, the notion of conductor of a character goes through, and
the functional equation of the Hecke $L$-function of a primitive
(``eigentlicher'' in Landau's paper) is
given in \cite[Theorem LVI]{Landau*18-2}.

\section{Some General Lemmas}

\subsection*{On Mellin transform}

\begin{lem}\label{smoothdecay}
  When $ \Re s =0$
  and $|\Im s | \ge 1$ or when $\Re s\ge 1/2$, we have
  \begin{equation*}
    |\check{w}(s)| ~\le~
    \frac{2^{\frac{n_\K}{2}+3}~\|w^{(n_\K+3)}\|_\infty }{(1 + |s|)^{n_\K+3}}.
  \end{equation*}
\end{lem}

\begin{proof}
  We set $A= n_{\K} + 3$ and $t=\Im s$ for typographical simplification. Integrating
  by parts $A$ times and noting that $w^{(m)}(1)= w^{(m)}(0)=0$ for
  $0 \le m\le A$, we get
  \begin{equation*}
    \left| \check{w}(s)=
      \frac{(-1)^A}{s(s+1)\cdots(s+A-1)}\int_{1/10}^1 w^{(A)}(u) ~u^{s+A-1}du \right|
   ~ \le~ 
    \frac{\|w^{(n_\K+3)}\|_\infty}{|s||s+1| \cdots |s+A-1|}.
  \end{equation*}
  In case $\Re s=0$, we note that $|t| \ge 1$ implies $1 + |t| \le 2|t|$.
  Furthermore we have $\sqrt{2}~|m + it | \ge |t|+1$ for all
  $m \ge 1$. This yields the constant $2^{\frac{n_\K}{2} + 2}$.

  In case $\Re s\ge 1/2$, we first notice that $3|s|\ge |s|+1$, and
  then we prove that $\sqrt{2}~|s+ m| \ge |s|+1$ for $m\ge1$ as before.
  To prove this inequality, set $\sigma=\Re s$. It is enough to prove
  that $t^2+2(\sigma+m)^2\ge \sigma^2+2|s|+1$. As
  $|s|\le |t|+|\sigma|$, it is enough to prove that
  $t^2-2|t|+1\ge -\sigma^2+2(1-2m)\sigma+2(1- m^2)$ which is obviously
  true.  This yields the constant $3\cdot2^{(n_\K+2)/2}$. We majorize
  the constant in both cases by $2^{\frac{n_\K}{2}+3}$ and this
  concludes this lemma.
 \end{proof}

%%%%%%%%%%%%%%
\begin{lem}
  \label{getM}
  For $\varepsilon\in(0,1/2]$, we define
  \begin{equation}
    \label{defM}
    M(w,\varepsilon)=\int_{-\infty}^\infty |\check{w}(it)|
    (1+|t|)^{\frac{1 +
        \varepsilon}{2}n_\K}dt.
  \end{equation}
  We have 
  $
    M(w,\varepsilon)
    \le
    2^{2 +\frac{\varepsilon n_\K}{2}}
    \bigr(\|w^{(n_\K+3)}\|_\infty 
    +
    10\cdot 2^{\frac{n_\K}{2}}\|w\|_1\bigl)$.
\end{lem}
%%%%%%%%%%%%%%

%%%%%%%%%%%
\begin{proof}
  Set $n=n_\K$.
  We split this integral according to whether $|t|\ge1$ or not. When
  $|t|\ge1$, Lemma~\ref{smoothdecay}  applies. When
  $|t|\le 1$, we simply use $|\check{w}(it)|\le 10 \|w\|_1$. This
  gives us
  \begin{equation*}
    M(w,\varepsilon)
    ~\le~
      2\frac{2^{\frac{n}{2}+3}\|w^{(n+3)}\|_\infty}
      {2^{\frac{1-\varepsilon}{2}n+2}(\frac{1-\varepsilon}{2}n+2)}
      +
      10\|w\|_12^{\frac{1+\varepsilon}{2}n+2}
    ~\le~
    2^{2+\frac{\varepsilon n}{2}}
    \bigr(\|w^{(n+3)}\|_\infty 
    +
    10\cdot 2^{\frac{n}{2}}\|w\|_1\bigl).
  \end{equation*}
  This completes the proof of the lemma.
\end{proof}
%%%%%%%%%%%

%%%%%%%%%%%%%%
\begin{lem}
  \label{getMstar}
  For $\varepsilon\in(0,1/2]$ and $r\in\{1,2\}$, we define
  \begin{equation}
    \label{defMstar}
    M^*(\varepsilon,r)=
    \int_{\frac{1+\varepsilon}{2}-i\infty}^{\frac{1+\varepsilon}{2}+i\infty} |\check{w}_0(s)|
    (1+|s|)^{\frac{1 +
        \varepsilon}{2r}n_\K}ds.
  \end{equation}
  We have 
  $M^*( \varepsilon,r) ~\le~ 12(57n_\K)^{n_\K + 3}$.
\end{lem}
%%%%%%%%%%%%%%

%%%%%%%%
\begin{proof}
  Set $n=n_\K$ and $\sigma=(1+\varepsilon)/2$. Lemma~\ref{smoothdecay}  applies and gives us
  \begin{equation*}
    M^*(\varepsilon,r)
    ~\le~
     2\cdot 2^{\frac{n}{2}+3}\|w_0^{(n+3)}\|_\infty
    \int_{0}^{ \infty} 
    \frac{dt}{(1+|\sigma + it|)^{\frac{(2r-1-\varepsilon)n}{2r}+ 3}}.
  \end{equation*}
  Furthermore
  $$
\int_{0}^{ \infty} 
 \frac{dt}{(1+|\sigma + it|)^{\frac{(2r-1-\varepsilon)n}{2r}+ 3}}
 ~\le~
  \int_{0}^{ \infty} 
 \frac{dt}{(1+t)^{\frac{(2r-1-\varepsilon)n}{2r}+ 3}}
 ~\le~
 1/2.
 $$
 Lemma~\ref{studyw0} leads to the bound $12(40\sqrt{2}n)^{n+ 3}$
 which is indeed not more than $12 (57 n)^{n + 3}$.
\end{proof}
%%%%%%%%

\subsection*{On the Dedekind zeta-function}

\begin{lem}\label{HR1}
 Let $0 < \varepsilon \le 1/2$.  In the strip $- \varepsilon \le \sigma\le 1 + \varepsilon$, 
 the Dedekind zeta-function $\zeta_{\K} (s)$ satisfies the inequality
 \begin{equation*}
 |\zeta_{\K}(s)| ~\le~ 3 \zeta(1+\varepsilon)^{n_{\K}}\biggl|\frac{s+1}{s-1}\biggr|
    \bigl(|d_{\K}|(1+|s|)^{n_{\K}}\bigr)^{\frac{1 + \varepsilon -\sigma}2}.
\end{equation*}
\end{lem}

%%%%%%%%%%%%
\begin{proof}
This is Theorem 4 of \cite{Rademacher*59} by Rademacher.
\end{proof}
%%%%%%%%%%%%
As a corollary, we deduce the following lemma. 
%%%%%%%%%%%%%
\begin{lem}
  \label{boundalphaK}
  We have
  $\displaystyle
    \frac{9\cdot 2^{n_\K}h_\K}{100 \sqrt{|d_\K|}}
    \le  \alpha_\K
    \le  6\left( \frac{2\pi^2}{5} \right)^{n_\K}|d_\K|^{1/4}$
    and
    $h_\K \le  67  (\pi^2/5 )^{n_\K} |d_\K|^{3/4}$.
\end{lem}
%%%%%%%%%%%%%
See the book \cite{Pohst-Zassenhaus*89} by Pohst and Zassenhaus for more on lower bounds for $\alpha_\K$.
%%%%%%%%
\begin{proof}
  Lemma~\ref{HR1} with the choices $\varepsilon=1/2$ and $s=1$ gives
  the upper bound.  As
  $2^{r_1}(2\pi)^{r_2}\ge 2^{r_1}2^{2r_2}\ge 2^{n_\K}$ and 
  the ratio of the regulator to $\mu_\K$ is bounded below absolutely
  by $0.09$ (see \cite{Friedman*89} by E.~Friedman), Eq.~\eqref{acf}
  provides us with the lower bound for $\alpha_\K$. Concerning the
  upper bound for $h_\K$, we use~\eqref{acf} again and this
  time, derive from it that $\alpha_\K\ge
  (9/100)\,2^{n_\K}h_\K/\sqrt{|d_\K|}$, from which the last bound
  follows immediately.
\end{proof}
%%%%%%%%

Next we deduce similar bounds for Hecke $L$-functions
corresponding to primitive characters $\chi$ of finite order.
For a Hecke character $\chi$ defined modulo $\q$,
the Hecke L-function associated to $\chi$ is defined as follows;
$$
L_{\q}(s,\chi) =  \sum_{\a \subseteq \O_{\K} \atop \a \ne 0} \frac{\chi(\a)}{\N(\a)^{s}},
$$
where $\Re s >1$.
We now state a result which bounds the growth of the
Hecke L-series using the Phragm{\'e}n-Lindel{\"o}f principle.

\begin{lem}[\cite{Rademacher*59}, Theorem 5]\label{HR2}
  Let $0 < \varepsilon \le 1/2$. In the strip
  $- \varepsilon \le \sigma\le 1 + \varepsilon$, the Hecke $L$-series
  associated with the primitive character $\chi$ of finite order and
  conductor~$\q$ satisfies the inequality
\begin{equation*}
 |L_{\q}(s,\chi)|
 ~\le~
  \zeta(1+ \varepsilon)^{n_{\K}}
\bigl(|d_{\K}|\mathfrak{N}(\q)(1+|s|)^{n_{\K}}\bigr)^{\frac{1+ \varepsilon -\sigma}2}.
\end{equation*}
\end{lem}
\begin{proof}
  This is a direct consequence of \cite[Theorem 5]{Rademacher*59} with
  $\eta= \varepsilon$, but one should be mindful of the notation,
  since Rademacher considers general Hecke Grossencharakteren, not
  necessarily of finite order. Things are rather clear when we inspect
  the gamma-factor given in~\cite[Bottom of page
  202]{Rademacher*59}. We have $a_p=0$ for every complex place,
  $v_p=0$ for every place and $a_p\in\{0,1\}$ for real places, $q$ of
  them taking the value 1.
\end{proof}

\subsection*{On rational primes}

%%%%%%%%%%%
\begin{lem}
  \label{takeiteasy}
  For $x \ge1$, we have
  $\displaystyle
    \sum_{\substack{p^k\le x, \atop k\ge2}}1
    ~\le~
    \frac{5\sqrt{x}}{4}.$
\end{lem}
%%%%%%%%%%%

%%%%%%%%%%%%%%%
\begin{proof}
  We first check this property by Pari/GP for $t\le 10^7$. To extend
  this result, let us denote by $S$ the sum to be bounded
  above. Inequality \cite[(3.32)]{Rosser-Schoenfeld*62} by Rosser and
  Schoenfeld tells us that
  \begin{equation}
    \label{rosserpsi}
    \sum_{p \le x} \log p ~\le~ 1.02\,x 
    \phantom{m}\text{for } x>0.
  \end{equation}
  We then readily check
  that
  \begin{align*}
    S  
    &\le~ \sum_{p\le \sqrt{x}}\frac{\log x}{\log p}
    ~\le~ 
    2(\log x)\int_2^{\sqrt{x}}\sum_{p\le u} \log p\frac{du}{u(\log
      u)^3}
    +(\log x)\sum_{p\le \sqrt{x}}\frac{\log p}{(\frac12\log x)^2}
    \\&
    \le~ 2.04(\log x)\int_2^{\sqrt{x}}\frac{du}{(\log u)^3}
    ~+~ 4.08~\frac{\sqrt{x}}{\log x}
    ~\le~ \sqrt{x}
  \end{align*}
  when $x \ge 10^7$. The lemma follows readily.
\end{proof}
%%%%%%%%%%%%%%%

%%%%%%%%%%%
\begin{lem}
  \label{takeiteasyb}
  For $x \ge100$, we have
  $\displaystyle
    \sum_{\substack{p\le x}}\frac{1}{p}
    ~\le~ 2\log\log x.
  $
\end{lem}
%%%%%%%%%%%
%%%%%%%%%%%
\begin{proof}
  This is readily checked by Pari/GP for $x \le 10^8$. We conclude the
  proof by appealing to \cite[Theorem 5]{Rosser-Schoenfeld*62} by Rosser
  and Schoenfeld.
\end{proof}
%%%%%%%%%%%

\subsection*{On the M\"{o}bius function}

For an ideal $\b$ of $\O_{\K}$, we define the M\"{o}bius
function as 
\begin{equation}
\label{defmoebius}
\mu(\b) = 
\begin{cases} 
1 &\text{ if  } \b = \O_{\K}, \\
(-1)^{r} &\text{ if  } \b =\P_1\cdots \P_r, ~\text{ where } \P_i \text{ are distinct prime ideals}, \\
0 &\text{ otherwise}. 
\end{cases}
\end{equation}
%%%%%%%%%%%%%%%%%
For a positive integer $R$ and for an ideal $\b$ of $\O_{\K}$, we define the truncated 
M\"{o}bius function $\mu_{R}$
in the following manner;
\begin{eqnarray*}
\mu_{R} (\b)
=
\begin{cases}
\mu(\b) & \text{ if } \omega(\b) \le R \\
0 & \text{ otherwise}.
\end{cases}
\end{eqnarray*}
Let $\psi_R(\b) = \sum_{\d \mid \b} \mu_{R}(\d)$. Applying M\"{o}bius
inversion formula and the fact that $\sum_{\d \mid \b}\mu(\d)=1$ if
and only if $\b = \O_{\K}$, we get 
$\mu_{R}(\b) = \sum_{\d \mid \b} \mu(\frac{\b}{\d}) \psi_R(\d)$. 
%%%%%%%%%%%%
In this context, we have the following Lemma. 
%%%%%%%%%%%%
\begin{lem}\label{boundpsi}
For any integral ideal $\d \neq \O_{\K}$, we have
$$
|\psi_{R} (\d)|
~\le~ 
{{\omega(\d) - 1} \choose R}.
$$
\end{lem}
\begin{proof}
Applying induction on $R$, we get
$$
\psi_{R} (\d)
~=~
\sum_{0 \le k \le R} (-1)^k {{\omega(\d) } \choose k}
~=~
(-1)^R {{\omega(\d) - 1} \choose R}
$$
(see H. Halberstam and H. Richert \cite[p. 46/47]{Halberstam-Richert*79}).
\end{proof}
\section{Degree one primes in ray class groups}
%%%%%%%%%%%%%%%%%%%%%%%%%%%%%%%%%%%%%%%%%%%%%%%%%%%%%%%%%%%%%%%%%%%%%%%%%%%%%%%%%%%
%%%%%%%%%%%%%%%%%%%%%%%%%%%%%%%%%%%%%%%%%%%%%%%%%%%%%%%%%%%%%%%%%%%%%%%%%%%%%%%%%%%

Let $\q$ be a modulus and let $T$ be a subgroup of index~2 of
$H_\q(\K)$. Let $\chi$ be the quadratic character whose kernel
is~$T$. We want to show that there exists a degree one prime $\P$ with
small norm such that $\chi(\P)=-1$ and another prime $\P'$ with small
norm such that $\chi(\P')=1$.

%%%%%%%%%%%%%%%%%%%%%%%%%%%%%%%%%%%%%%%%%%%%%%%%%%%
\subsection*{$L$-series for degree one primes}
%%%%%%%%%%%%%%%%%%%%%%%%%%%%%%%%%%%%%%%%%%%%%%%%%%%

For a Hecke character of finite order modulo $\q$, we define
\begin{equation}\label{defFchi}
F(s, \chi)  ~=~ \prodflat_{\substack{ \P \nmid \q}}~ \frac{1}{1 - \chi(\P)\N(\P)^{-s}}.
\end{equation}
%%%%%%%%%%%
\begin{lem}\label{boundFchi}
When $\Re s = (1+ \varepsilon)/2$ for some $0< \varepsilon  \le 1/2$ and 
$\chi$ is a Hecke character of finite order modulo~$\q$, we have
$$
|F(s,\chi)| ~\le~ \zeta(1+\varepsilon)^{\frac{3n_{\K}}{2}} (|d_{\K}|\N(\q))^{\frac{1+ \varepsilon}{4}} 
\theta(\q) (1+|s|)^{\frac{(1+\varepsilon)n_{\K}}{4}},
$$
where $\theta(\q)=\prod_{\P|\q}\frac{\sqrt{\N\P}}{\sqrt{\N\P}-1}$.  When
$\chi = \chi_{0, \q}$ is the trivial character modulo $\q$, we have
$$
|F(s,\chi_{0, \q})| ~\le~  27 \zeta(1+\varepsilon)^{\frac{3n_{\K}}{2}} |d_{\K}|^{\frac{1+ \varepsilon}{4}} 
\theta(\q) (1+|s|)^{\frac{(1+\varepsilon)n_{\K}}{4}}.
$$
\end{lem}
%%%%%%%%%
%%%%%%%%%%%%%%%%
\begin{proof}
To find an upper bound for $F(s,\chi)$, we first reduce it to Hecke
$L$-series using the following product $F(s,\chi) = L_{\q}(s,\chi)J(s,\chi)$,
where
\begin{equation}\label{defJ}
  J(s,\chi) ~=~ \prodsharp_{\substack{ \P \nmid \q}} ~\bigl(1-\chi(\P)\N(\P)^{-s}\bigr).
\end{equation}
We readily find that, when $\Re s = (1+ \varepsilon)/2$, we have
\begin{equation}
  \label{eq:7}
  |J(s,\chi)| ~\le~ \prod_{p}\bigl(1+p^{-2\sigma}\bigr)^{\frac{n_{\K}}{2}}
  ~\le~ 
  \zeta(1+\varepsilon)^{ \frac{n_{\K}}{2} }.
\end{equation}
The next step is to reduce $L_{\q}(s,\chi)$ to
$L_{\mathfrak{f}}(s,\chi^*)$, where $\chi^*$ is the primitive character, say 
modulo~$\mathfrak{f}$, inducing~$\chi$. We directly see that
$$
L_{\q}(s,\chi) = L_{\mathfrak{f}}(s,\chi^{*})\prod_{\P \mid \q \atop{ \P  \nmid \mathfrak{f} }} 
(1 - \chi(\P)\N(\P)^{-s}).
$$
Therefore, applying Lemma~\ref{HR2} when $\chi^*$ is not equal to the
constant character~1, and
Lemma~\ref{HR1} otherwise, together with the bound in \eqref{eq:7}, we get
the desired result.
\end{proof}
%%%%%%%%%%%

%%%%%%%%%%%%%%%%%%%%%%%%%%%%%%%%%%%%%%%%%%%%%%%%%%%%%%%%%%%%%%%%%%%%%%%%%%%%%%%%%%% 
\subsection*{Products of degree one primes in ray classes modulo $\q$,
sieve approach}
%%%%%%%%%%%%%%%%%%%%%%%%%%%%%%%%%%%%%%%%%%%%%%%%%%%%%%%%%%%%%%%%%%%%%%%%%%%%%%%%%%%

%%%%%%%%%%%
\begin{lem}
  \label{primeboundini}
  Let $\mathfrak{b}$ be an integral ideal co-prime to $\q$.
  When $F_1(\q) = 2^{r_1} h_\K \N\q$ and
  \begin{equation*}
    X 
  ~\ge~
  \log(3F(\q))^{n_\K^2}\, B(\K)F_1(\q)\N\q \, \log\log(B(\K)F(\q)\N\q)^2,
  \quad
  B(\K) = (n_\K^{50n_\K^2} E(\K)\sqrt{|d_\K|})^{n_\K},
 \end{equation*}
  we have
\begin{equation*}
  \sum_{\substack{\P  \\ \N\P\ge 40^{n_\K} \atop { \deg \P\ge 2}}  }
  \sum_{\substack{\P|\mathfrak{a} \subseteq \O_{\K}
    \\ [\mathfrak{a}] 
      = [\mathfrak{b}] \\ \N\mathfrak{a}\le X}}
  1
~~\le~
\frac{\alpha_\K\phi(\q)X}{2^{{n_\K} } \N\q ~|H_\q(\K)|}.
\end{equation*}
\end{lem}
%%%%%%%%%%%
We can remove the term $\log\log(B(\K)F(\q)\N\q)^2$ when $n_\K \ge 3$.

%%%%%%%%%%%
\begin{proof}
  We shorten $n_\K$ by $n$ and set $D=40^n$.
  On calling $S$ the sum on the left hand side, Theorem~\ref{asymfinal}
  gives us the upper bound, with $Y=\alpha_{\K} \phi(\q) X /(\N\q h_{\K,\q})$,
  \begin{equation*}
    S 
    ~\le~
    \sum_{\substack{\P \\ \N\P\ge D,\\ \deg \P\ge2}}\frac{Y}{\N\P}
    ~+~
    E(\K)
    F(\q)^{\frac{1}{n}}\log(3F(\q))^{n}
    \sum_{\substack{\P \\ X\ge\N\P\ge D,\\ \deg
        \P\ge2}}\left(\frac{X}{\N\q\N\P}\right)^{1-\frac{1}{n}}
    ~+~
    n^{8n}\frac{R_\K}{|\mu_\K|} F(\q)\sum_{\substack{\P \\ \N\P\le X,\\ \deg
        \P\ge2}}1.
  \end{equation*}
  We now have to study each of these three sums. We readily bound
  above the first one by
  \begin{equation*}
    Y\sum_{\substack{p^k\ge D,\\ k\ge2}}\frac{n/2}{p^k}
    ~\le~ \frac{nY}{2}\int_{D}^{\infty} (\sum_{\substack{D\le p^k\le t,\\ k\ge2}}1 )~\frac{dt}{t^2}
    ~\le~ \frac{5nY}{8}\int_{D}^{\infty}\frac{dt}{t^{3/2}}
    ~\le~ \frac{5nY}{4\sqrt{D}}
  \end{equation*}
  by Lemma~\ref{takeiteasy}. We notice that $(5/4)40^{-n/2}\le
  1/2^{2n+1}$ when $n\ge2$.
  % f(n)=5/4*2^(n+1)/10^(n/2)
  The same Lemma~\ref{takeiteasy} yields
  the bound $\tfrac58 n^{8n+1}R_\K F(\q)\sqrt{X}/|\mu_\K|$ for the third term. We
  further find that
  \begin{align*}
    %\notag
    \tfrac58 n^{8n+1}\frac{R_\K}{|\mu_\K|} F(\q)\sqrt{X}
    &=~\frac{5n^{8n+1}}{8} \frac{R_\K \N\q \, h_{\K,\q} F(\q) }{\alpha_\K |\mu_\K| \phi(\q)} 
   \frac{Y\sqrt{X}}{X}
    ~=~  \frac{5n^{8n+1}}{8} \frac{2^{r_1} h_{\K}R_\K}{\alpha_\K|\mu_\K|}
      \frac{\N\q \,Y}{\sqrt{X}}
    \\&\le~%\label{err1}
  \frac{5n^{8n+1}}{8}  (2\pi)^{-r_2}\sqrt{|d_\K|}
~ \frac{\N\q \, Y}{\sqrt{X}}
  ~\le~  n^{9n}\sqrt{|d_\K|} ~\frac{\N\q \, Y}{\sqrt{X}} 
  \end{align*}
 by applying the definition of $F(\q)$ from Theorem~\ref{asymfinal}
 and the expression of $\alpha_\K$ in terms of the invariants of the
  field mentioned in~\eqref{acf}. Let us now examine the second term
  above, a task for which we distinguish between the case when $n=2$
  and $n\ge3$. In the latter situation, we find that
  \begin{equation}
    \label{compl}
    \sum_{\substack{\P \\ X\ge\N\P\ge D,\\ \deg
        \P\ge2}}\frac{1}{(\N\P)^{1-\frac{1}{n}}}
    ~\le~ \frac{n}{2}\sum_{p\ge2}\frac{1}{p^{2(1-1/3)}}
    ~\le~ n
  \end{equation}
  by appealing to Pari/GP.
  In this case, we thus find that
  \begin{align*}
    S/Y
    &\le~ \frac{n}{2^{2n+1}}
    +
    E(\K)F(\q)^{\frac{1}{n}}\log(3F(\q))^{n}
    \frac{nh_{\K,\q}}{\alpha_\K\phi(\q)}\left(\frac{\N\q}{X}\right)^{{1}/{n}}
    ~+~
   n^{9n}\sqrt{|d_\K|} \frac{\N\q}{\sqrt{X}}
    \\&\le~
     \frac{1}{2^{n+1}}
    +
    E(\K)\log(3F(\q))^{n}
    \frac{n2^{r_1}h_{\K}}{\alpha_\K}\left(\frac{F(\q)\N\q}{X}\right)^{{1}/{n}}
    ~+~
  n^{9n}\sqrt{|d_\K|}  \frac{\N\q}{\sqrt{X}}
    \\
    &\le~ \frac{1}{2^{n+1}}
    +
     \frac{100n}{9}E(\K)\log(3F(\q))^{n}
  \sqrt{|d_\K|}\left(\frac{F(\q)\N\q}{X}\right)^{{1}/{n}}
    ~+~
   n^{9n}\sqrt{|d_\K|} \frac{\N\q}{\sqrt{X}}
  \end{align*}
  by Lemma~\ref{classnumber} and~\ref{boundalphaK}. Since we have assumed
  that
  \begin{equation}
    \label{eq:2}
    X~\ge~ B(\K)F_1(\q)\N\q \,\log(3F(\q))^{n^2} ,
    \quad
    B(\K)  = (n^{50n^2} E(\K)\sqrt{|d_\K|})^{n}
  \end{equation}
  we find
  that $S\le Y/2^{n}$ when $n\ge3$. When $n=2$, we only have to replace
  the upper bound $n$ in \eqref{compl} by $n\log\log X=2\log\log X$ by
  Lemma~\ref{takeiteasyb}. Proceeding as above we reach the inequality
  \begin{equation*}
    S/Y
    ~\le~ \frac14
    ~+~
   \frac{200}{9} E(\K)\log(3F(\q))^2
    \log\log X \sqrt{|d_\K|}\left(\frac{F(\q)\N\q}{X}\right)^{{1}/{2}}
    ~+~
  2^{18}  \frac{\N\q}{ \sqrt{X}}\sqrt{|d_\K|}.
  \end{equation*}
  Our hypothesis on $X$ reads
  \begin{equation*}
    X ~\ge~ \log(3F(\q))^{4}\,B(\K)F_1(\q)\N\q\,\log\log(B(\K)F(\q)\N\q)^2.
  \end{equation*}
  We just need to notice that
  $\log(3F(\q))^{4}\log\log(B(\K)F(\q)\N\q)^2\le
  (B(\K)F(\q)\N\q)^6$, so that
  \begin{align*}
    \log\log X
    &\le \log\log((B(\K)F(\q)\N\q)^7)
    \\&\le \biggl(1+ \frac{\log 7}{\log\log(B(\K))}\biggr)\log\log(B(\K)F(\q)\N\q)
    \le 2\log\log(B(\K)F(\q)\N\q).
  \end{align*}
  This is enough to conclude.
\end{proof}
%%%%%%%%%%%

%%%%%%%%%%%
\begin{lem}
  \label{primeboundinib}
  Let $\mathfrak{b}$ be an integral ideal co-prime to $\q$.
  When $F_1(\q) = 2^{r_1} h_\K \N\q$ and
  \begin{equation*}
    X ~\ge~ 
  \log(3F(\q))^{n_\K^2}\, B(\K)F_1(\q)\N\q \, \log\log(B(\K)F(\q)\N\q)^2,\quad
  B(\K)=(n_\K^{50n_\K^2} E(\K)\sqrt{|d_\K|})^{n_\K},
  \end{equation*}
  we have
\begin{equation*}
  \sumflat_{\substack{\mathfrak{a} \subseteq \O_{\K}
    \\ [\mathfrak{a}] 
      = [\mathfrak{b}] \\ \N\mathfrak{a}\le X}}~
  1
~\ge~
(1/4)\frac{\alpha_\K\phi(\q)X}{ \N\q |H_\q(\K)|}.
\end{equation*}
\end{lem}
%%%%%%%%%%%

%%%%%%%%%%%%%%%%%
\begin{proof}
  We denote $n_\K$ by $n$ and set $D=40^n$. Let $M$ be the product of
  the prime ideals in $ \O_{\K}$ of degree greater than or equal to $2$, 
  co-prime to $\q$ and of norm at most $D$.
Denoting the sum to be evaluated by $S$, a simple combinatorial
  argument together with Lemma~\ref{primeboundini} gives us
  \begin{equation*}
    S
    ~\ge~
    \sum_{\substack{\mathfrak{a} \subseteq \O_{\K}
        \\ [\mathfrak{a}] = [\mathfrak{b}] 
        \\ \N\mathfrak{a}\le X\\ (\mathfrak{a},M)=1}}
    1
    ~- ~\frac{\alpha_\K\phi(\q)X}{2^{n} \N\q |H_\q(\K)|}
    ~=~
    S(M) ~-~  \frac{\alpha_\K\phi(\q)X}{2^{n} \N\q |H_\q(\K)|}
  \end{equation*}
  say. We detect the coprimality condition with an elementary instance
  of the Brun sieve. We select an odd positive integer $R$ and write
  \begin{equation*}
    1_{(\mathfrak{a},M)=1}=\sum_{\mathfrak{d}|(M,\mathfrak{a})}\mu(\mathfrak{d})
    \ge \sum_{\substack{\mathfrak{d}|(M,\mathfrak{a}),\\
        \omega(\mathfrak{d})\le R}}\mu(\mathfrak{d}).
  \end{equation*}
  We deduce
  \begin{equation*}
    S(M) ~\ge~
    \sum_{\substack{\mathfrak{d}|M,\\
        \omega(\mathfrak{d})\le R}}\mu(\mathfrak{d})
    \sum_{\substack{\mathfrak{a}' \subseteq \O_{\K}
        \\ [\mathfrak{d}\mathfrak{a}'] 
        = [\mathfrak{b}] \\ \N\mathfrak{a}'\le \frac{ X}{ \N\mathfrak{d}} }}1.
  \end{equation*}
  Theorem~\ref{asymfinal}
  gives us the lower bound, with $Y=\alpha_{\K} X\phi(\q)/(\N\q \,
  h_{\K,\q})$,
  \begin{equation*}
    S(M)
    ~\ge~
    Y\sum_{\substack{\mathfrak{d}|M,\\
        \omega(\mathfrak{d})\le R}}\frac{\mu(\mathfrak{d})}{\N\mathfrak{d}}
    -
    E(\K)
    F(\q)^{\frac{1}{n}}\log(3F(\q))^{n}
    \sum_{\substack{\mathfrak{d}|M,\\
        \omega(\mathfrak{d})\le R}}
    \left(\frac{X}{\N\q\N\mathfrak{d}}\right)^{1-\frac{1}{n}}
    -~
 n^{8n}\frac{R_\K}{|\mu_\K|} F(\q) \,  
 \sum_{\substack{\mathfrak{d}|M,\\  \omega(\mathfrak{d})\le R}}1.
  \end{equation*}
 Concerning the main term, we can write
 \begin{eqnarray}\label{New}
 \sum_{\substack{\mathfrak{d}|M,\\  \omega(\mathfrak{d})\le R}}
\frac{\mu(\mathfrak{d})}{\N\mathfrak{d}}
~=~ 
\sum_{\d \mid M} \frac{\mu_{R} (\d)}{\N\d}  
~=~
\sum_{\d \mid M} \frac{1}{\N\d} \sum_{\b \mid \d} \mu \left(\frac{\d}{\b}\right) \psi_R(\b) 
&=&
\sum_{ \b \mid M} \frac{\psi_R(\b)}{\N\b} \sum_{\c \mid \frac{M}{\b}} \frac{\mu(\c)}{\N\c} \nonumber \\
&=&
\prod_{\P \mid M} \left( 1- \frac{1}{\N\P} \right) \sum_{ \b \mid M} \frac{\psi_R(\b)}{\phi(\b)}.
\end{eqnarray}
Applying \lemref{boundpsi}, we have
$$
\sum_{\b \mid M \atop \N\b > 1}  \frac{\psi_R(\b)}{\phi_{\K}(\b)} 
~\le~
\sum_{\b \mid M \atop \N\b > 1} {{\omega(\b) - 1} \choose R} \frac{1}{\phi(\b)}
~\le~
\sum_{m=R+1}^{\pi_{\K}(40^{n})} {{m - 1} \choose R}
\sum_{\b \mid M \atop {\N\b > 1 \atop \omega(\b)=m}} \frac{1}{\phi(\b)},
$$
where $\pi_{\K}(x)$ be the number of prime ideals in $\O_{\K}$ with norm less
than or equal to $x$. The last quantity is equal to
\begin{eqnarray*}
\sum_{m=R+1}^{\pi_{\K}(40^{n})} {{m - 1} \choose R} \left( \sum_{\P \mid M} 
\frac{1}{\phi(\P)} \right)^m  \frac{1}{m!}
&\le& 
\frac{1}{R!} \sum_{m=R+1}^{\pi_{\K}(40^{n})} \frac{1}{(m-R)!}  
\left( \sum_{\P \mid M} \frac{1}{\N\P-1} \right)^m \\
&\le&
\frac{1}{R!} \sum_{m=R+1}^{\pi_{\K}(40^{n})} \frac{1}{(m-R)!}  
\left( \sum_{p < 40^{n}}
 \frac{n}{2(p^2-1)} \right)^m \\
&\le& 
\frac{1}{R!} \sum_{m=R+1}^{\pi_{\K}(40^{n})} \frac{1}{(m-R)!}  
\left( \sum_{p < 40^{n}} \frac{{n}}{p^2} \right)^m.
\end{eqnarray*}
We have, however, that
$$
\sum_{p < 40^{n}} \frac{1}{p^2} ~\le~ \int_{1}^{40^{n}} \frac{dt}{t^2} ~\le~ 1.
$$
Therefore, we get
\begin{eqnarray*}
\sum_{\b \mid M \atop \N\b > 1}  \frac{\psi_R(\b)}{\phi(\b)}
& \le &
\frac{1}{R!} \sum_{m=R+1}^{\pi_{\K}(40^{n})} \frac{n^m}{(m-R)!}  
~\le~ 
\frac{n^R }{R!} \sum_{m=R+1}^{\pi_{\K}(40^{n})} \frac{ {n}^{m-R}}{(m-R)!} 
~\le~ 
\frac{ n^R e^{n} }{R!}  
~\le~ 
\frac{e^{n + R\log n }}{R!}.
\end{eqnarray*}
We know that $ R! \ge (\frac{R}{e})^R$ and this implies that
$$
\sum_{\b \mid M \atop \N\b > 1}  \frac{\psi_R(\b)}{\phi(\b)}
~\le~ 
\frac{e^{n + R\log n }}{R^R}e^R
~\le~
\exp(n + R - R\log R +R \log {n}).
$$
 We
  select
  \begin{equation}
    \label{defR}
    R=2[5n\log n]+1
  \end{equation}
  so that
\begin{equation}\label{NEWF}
\exp \left(n - R \log \frac{R}{en} \right) 
~\le~ 
\exp\left(n  - (9n\log n ) \log \frac{9\log n }{e} \right)
~\le~
\exp\left(- 4 n \right)
~\le~ 
50^{-n}.
\end{equation}
Combining \eqref{New} and \eqref{NEWF}, we get 
$$
 \sum_{\substack{\mathfrak{d}|M,\\
        \omega(\mathfrak{d})\le R}}
        \frac{\mu(\mathfrak{d})}{\N\mathfrak{d}} 
~\ge~ 
\prod_{\P \mid M} \left( 1- \frac{1}{\N\P} \right) (1 - 50^{-n})
~\ge~
 \left(\frac{6}{\pi^2} \right)^{\frac{n}{2}} \left(1 - 10^{-n} \right),
$$
since all the prime ideals dividing $M$ have degree at least $2$.
This takes care of the main term. Concerning the error term for
  $S(M)$, we first notice that the number of rational primes less than $D$ is
  at most $3D/(2\log D)$ (see \cite[(3.6)]{Rosser-Schoenfeld*62}). This implies 
  that the number of prime ideals in
  $\O_{\K}$ of norms at most $D$ and of degree $\ge2$ is at most
  $[n/2]3D/(2\log D)$ which is not more that $2D/3$. Whence
  \begin{equation*}
    \sum_{\substack{\mathfrak{d}|M,\\
        \omega(\mathfrak{d})\le R}}1
    ~\le~ (2D/3)^R.
  \end{equation*}
  We have thus reached the inequality
  \begin{equation*}
    S(M)
    ~\ge~
    Y\bigl(1-{10^{-n}}\bigr)(6/\pi^2)^{n/2}
    ~-~
    (2D/3)^R\biggl(E(\K)
    F(\q)^{\frac{1}{n}}\log(3F(\q))^{n}
    \left(\frac{X}{\N\q}\right)^{1-\frac{1}{n}}
    ~+~
    n^{8n}\frac{R_\K}{|\mu_\K|} F(\q)\biggr).
  \end{equation*}
  We  get
  \begin{align*}
    \frac{1}{Y}(2D/3)^RE(\K)
    F(\q)^{\frac{1}{n}}\log(3F(\q))^{n}
    \left(\frac{X}{\N\q}\right)^{1-\frac{1}{n}}
    &\le~
      40^{12n^2\log n}\frac{ h_{\K,\q}}{\alpha_{\K}
      \phi(\q)}\left(\frac{F(\q)\N\q}{X}\right)^{\frac{1}{n}}
      E(\K)
      \log(3F(\q))^{n}
    \\&\le~
    40^{12n^2 \log n}\frac{2^{r_1}h_\K}{\alpha_{\K}}\left(\frac{F(\q)\N\q}{X}\right)^{\frac{1}{n}}
      E(\K)
    \log(3F(\q))^{n}
    \\&\le~ 
    40^{12n^2 \log n}\frac{100}{9} \sqrt{|d_\K|}\left(\frac{F(\q)\N\q}{X}\right)^{\frac{1}{n}}
      E(\K)
  \end{align*}
  by Lemma~\ref{classnumber} and~\ref{boundalphaK}. Our lower bound
  on~$X$ ensures that this upper bound is $\le 1/10^n$.
  The second error term is treater similarly. We write
  \begin{align*}
    \frac{D^Rn^{8n}}{Y}\frac{R_\K}{|\mu_\K|} F(\q)
    &\le~
      \frac{40^{12n^2\log n}n^{8n}\N\q h_{\K,\q}}{\alpha_{\K}
      X\phi(\q)}\frac{R_\K}{|\mu_\K|} F(\q)
    \\&=~
    \frac{40^{12n^2\log n}n^{8n}\N\q h_{\K,\q}}{
      Xh_\K\phi(\q)2^{r_1}(2\pi)^{r_2}}
    \sqrt{|d_\K|}\frac{2^{r_1}\phi(\q)h_\K}{h_{\K,\q}}
    ~\le~
    40^{12n^2\log n}n^{8n}\frac{\N\q}{X}
    \sqrt{|d_\K|}
  \end{align*}
  by using \eqref{acf} and the definition of $F(\q)$. The lower bound
  for $X$ again implies that this quantity is at most $1/10^n$.
Combining all our estimates, we find that
  \begin{equation*}
    S/Y\ge
    \bigl(1-{10^{-n}}\bigr)(6/\pi^2)^{n/2}
    -\frac{1}{10^n}
    -\frac{1}{10^n}
    -\frac{1}{2^n}
    \ge
    \tfrac14(6/\pi^2)^{n/2}
  \end{equation*}
  as required. This concludes the proof.
\end{proof}
%%%%%%%%%%%%%%%%%

%%%%%%%%%%%%%%%%%%%%%%%%%%%%%%%%%%%%%%%%%%%%%%%%%%%%%%%%%%%%%%%%%%%%%%%%%%%%%%%%%%% 
\subsection*{Products of degree one primes in ray classes modulo $\q$,
analytic approach}
%%%%%%%%%%%%%%%%%%%%%%%%%%%%%%%%%%%%%%%%%%%%%%%%%%%%%%%%%%%%%%%%%%%%%%%%%%%%%%%%%%%

%%%%%%%%%%%
\begin{lem}
  \label{primebound}
  Let $\mathfrak{b}$ be an integral ideal co-prime to $\q$.
  When $X \ge  10^{25n_\K} n_{\K}^{7n_{\K}}
    |d_\K|^{4}  \N\q^3$, we have
\begin{equation*}
\sumflat_{\mathfrak{a} \subseteq \O_{\K} \atop{ [\mathfrak{a}] 
= [\mathfrak{b}]}} w_0\left(\frac{\N(\mathfrak{a})}{X}\right) 
~\ge~
\frac{\alpha_\K\phi(\q)\check{w}_0(1)X}{2(1.3)^{n_\K} \N\q ~|H_\q(\K)|}.
\end{equation*}
\end{lem}
%%%%%%%%%%%

%%%%%%%%%%%
\begin{proof}
On calling $S$ the sum on the left hand side, the orthogonality of characters 
readily gives us that (recall~\eqref{defFchi})
\begin{align*}
S 
&=~
\frac{1}{|H_{\q}(\K)|} \sum_{\chi \in \hat{H}_{\q}(\K)} \overline{\chi}(\mathfrak{b})
 \sumflat_{(\mathfrak{a}, \q)=\O_{\K}} 
     w_0\left(\frac{\N(\mathfrak{a})}{X}\right) \chi(\mathfrak{a}),
  \\
& =~ 
\frac{1}{|H_{\q}(\K)|} \sum_{\chi \in\hat{H}_{\q}(\K)} 
\frac{\overline{\chi}(\b)}{2\pi i} \int_{2-i\infty}^{2+i\infty} F(s,\chi)\check{w}_0(s) X^{s}ds.
\end{align*}
On using Lemma~\ref{boundFchi}, we find that
\begin{align*}
S 
& =~  
     \frac{\alpha_{\K}\phi(\q)J(1, \chi_{0, \q})}{\N\q}
     \frac{\check{w}_0(1)X}{|H_{\q}(\K)| } 
~+~
\frac{1}{|H_{\q}(\K)|}\sum_{\chi \in\hat{H}_{\q}(\K)} \frac{\overline{\chi}(\b)}{2\pi i} 
\int_{\frac{1+\varepsilon}{2}-i\infty}^{\frac{1+\varepsilon}{2}+
     i\infty} F(s,\chi)\check{w}_0(s) X^{s}ds
  \\
& =~  
\frac{\alpha_{\K}\phi(\q)J(1, \chi_{0, \q})}{\N\q}
       \frac{\check{w}_0(1)X}{|H_{\q}(\K)| } +  
\rO^*\left( 5 \zeta(1+\varepsilon)^{\frac{3n_{\K}}{2}} 
(|d_{\K}|\N\q)^{\frac{1+\varepsilon}{4}}\theta(\q) X^{\frac{1+\varepsilon}{2}} 
M^*(\varepsilon,2)
     \right)
\end{align*}
where $M^*$ is as defined in Lemma~\ref{getMstar}.
Proceeding as in~\eqref{eq:7}, we find that 
\begin{equation}\label{lb-j}
J(1, \chi_{0, \q})\ge \zeta(2)^{-n_\K/2}
\end{equation}
while Lemma~\ref{boundalphaK} provides us with a
lower bound for $\alpha_\K$ and, combined with~\eqref{eq:4}, an upper bound 
for $|H_\q(\K)|$. This together with Lemma~\ref{studyw0} tells us that the main term above
is at least, with the notation $n=n_\K$,
\begin{equation*}
  \frac{9X}{33500 \sqrt{\zeta(2)^{n}}(\pi^2/5)^{n} \sqrt{n}|d_\K|^{5/4}\N\q }
  ~\ge~
  \frac{X}{3723(\pi^2/5)^{3n/2} \sqrt{n}|d_\K|^{5/4}\N\q}.
\end{equation*}
This is larger than twice the above error term provided we have
\begin{equation}
  \label{ineqstep}
  X^{\frac{1-\varepsilon}{2}}
  ~\ge~
   83 \cdot 10^{9} n^{\frac{7}{2}} (159n)^{n}
    \zeta(1+\varepsilon)^{\frac{3n}{2}} 
    |d_\K|^{\frac{6+\varepsilon}{4}}\theta(\q) \N\q^{\frac{5+\varepsilon}{4}}.
  \end{equation}
  We select $\varepsilon=1/10$. The previous inequality is implied by
  \begin{equation*}
    X 
    ~\ge~
   19 \cdot 10^{23} n^{\frac{70}{9}} (5476n)^{\frac{20n}{9}}
    |d_\K|^{\frac{61}{18}}\theta(\q)^{\frac{20}{9}} \N\q^{\frac{17}{6}}.
  \end{equation*}
  Now $(\frac{\sqrt{p}}{\sqrt{p}-1})^{20/9}\le p^{1/6}$ when $p >19$, while
$
  \prod_{p\le
    19}\theta(p)^{20/9}/{p^{1/6}}
  \le 1200$.
As a conclusion, we derive that $\theta(\q)^{20/9}\le 1200^{n}\N\q^{1/6}$.
Some numerical computations end the proof of our lemma.
\end{proof}
%%%%%%%%%%%

%%%%%%%%%%%%%
\begin{proof}[Proof of Theorems~\ref{degreeoneprime} and \ref{degreeoneprimebis}]
  Theorem~\ref{degreeoneprime} follows as an easy consequence of
  Lemma~\ref{primebound}, and similarly, Theorem~\ref{degreeoneprimebis} follows
  from Lemma~\ref{primeboundinib}.
\end{proof}
%%%%%%%%%%%%%

\subsection*{Degree one primes in quadratic ray subgroups modulo $\q$}

%%%%%%%%%%%%%%%%%%%%%%%%%%%%%%%%%%%
\begin{lem}\label{boundLonechiquad}
Let $\chi$ be quadratic character on $H_\q(\K)$. We have
\begin{equation*}
L_{\q}(1,\chi)  
~\ge~  
\frac{9\cdot 2^{2n_\K}}{100 \alpha_{\K} |d_{\K}| \sqrt{\N(\q)}} 
\phantom{mm}\text{and}\phantom{mm}
F(1,\chi)  
~\ge~  
\frac{9\cdot 2^{2n_\K}}{100 \alpha_{\K} |d_{\K}| \sqrt{\zeta(2)^{n_{\K}}\N(\q)}} 
\end{equation*}
where $d_{\K}$ is the discriminant of $\K$,
$\alpha_{\K}$ is the residue of the Dedekind zeta function
at $s=1$. 
\end{lem}
%%%%%%%%%%%%%%%%%%%%%%%%%%%%%%%%
This lemma improves on \cite[Lemma 2]{Hinz-Lodemann*94} of Hinz and
Lodemann in two ways: the dependence in the base field is explicit
and the lower bound is in $1/\sqrt{\N\q}$ instead of $1/[\sqrt{\N\q}(\log\N\q)^2]$.

\begin{proof}
We note that the product $L_{\f}(s,\chi) \zeta_{\K}(s)$, where $\f | \q$
is the Dedekind zeta function of a quadratic extension $\M$
of $\K$. By \eqref{acf}, we have 
$$
L_{\f}(1,\chi)\alpha_{\K} = \alpha_{\M} = \frac{2^{r_1}(2\pi)^{r_2}h_{\M} 
R_{\M}}{|\mu_{\M}| \sqrt{|d_{\M}|}},
$$
where $r_1$ and $2r_2$ are the number of real and complex embeddings of $\M$. Since this field is of
degree $2n_\K$ over $\mathbb{Q}$, Lemma~\ref{boundalphaK} gives us
\begin{equation}\label{lbd-res}
  \alpha_{\M} ~\ge~ \frac{9\cdot 2^{n_\M}}{100 \sqrt{|d_{\M}|}}.
\end{equation}
Finally, since the extension $\M/\K$ has conductor with finite part $\f$
and is a quadratic extension, by the conductor discriminant
formula (\cite[Chapter VII, Point (11.9)]{Neukirch*99} by Neukirch) the relative 
discriminant of $\M/\K$ is also $\f$.
However we also have that (\cite[Chapter III, Corollary (2.10)]{Neukirch*99})
$\sqrt{|d_{\M}|} = |d_{\K}| \sqrt{\N(\f)}$. This gives that
$$
|L_{\q}(1,\chi)| ~\ge~  |L_{\f}(1,\chi) | \prod_{\P| \q \atop{\P\nmid \f} }\left(1 - \frac{1}{\N\P}\right).
$$
Since 
$$
\frac{\sqrt{\N\q}}{\sqrt{\N\f}} \prod_{\P| \q \atop{\P\nmid \f} }\left(1 - \frac{1}{\N\P}\right)
~\ge~
\frac{\sqrt{\N\q} \prod_{\P| \q} \left(1 - \frac{1}{\N\P}\right)}{\sqrt{\N\f}\prod_{\P| \f} \left(1 - \frac{1}{\N\P}\right)} 
~\ge~ 1,
$$
we get the desired bound.
To extend it to $F(1,\chi)$, we notice that
\begin{equation}\label{J-bound}
\left| \frac{L_{\q}(1,\chi)}{F(1,\chi)} \right|
    ~\le~
     \prodsharp_{\substack{\P \nmid \q}}
    ~\left(1- \frac{1}{\N\P} \right)^{-1}
    ~\le~ 
    \zeta(2)^{\frac{n_{\K}}{2}}. 
\end{equation}
\end{proof}

%%%%%%%%%%%
\begin{lem}\label{simAxer}
Let $\chi$ be a quadratic character on $H_\q(\K)$. We have 
$$
\sumflat_{\substack{\mathfrak{a} \subseteq \O_K \\ 
(\mathfrak{a},\q)=\O_K }}
 ({\bf 1}\star\chi)(\mathfrak{a})~w_0\left(\frac{\N(\mathfrak{a})}{ X}\right)
~ >~
 \frac{X}{27 
    \sqrt{\N\q}|d_\K|^2}\frac{\phi(\q)}{\N\q}
 $$
provided that $X\ge 8 \cdot (10^{31} n_\K^7)^{n_\K}  |d_\K|^4 \N\q^2$.
\end{lem}
%%%%%%%%%

%%%%%%%%%%%%
\begin{proof}
Let us denote the sum on the left hand side by $S_1(w_0)$ and the principal Hecke character 
modulo~$\q$ by $\chi_{0,\q}$.  By mimicking the proof of 
Lemma~\ref{primebound}, we find that
 \begin{align*}
   S_1(w_0)
   &=
     \frac{1}{2\pi i}\int_{2-i\infty}^{2+i\infty}F(s,\chi_{0,\q})F(s,\chi)\check{w}_0(s)X^s
     ds
   \\
   &=
     \alpha_{\K, \q} F(1,\chi)\check{w}_0(1) X + \frac{1}{2\pi i}
     \int_{\frac{1+\varepsilon}{2} - i\infty}^{\frac{1+\varepsilon}{2} + i\infty}
     F(s,\chi_{0,\q})F(s,\chi)\check{w}_0(s)X^s ds
   \\
   &=
     \alpha_{\K, \q} F(1,\chi)\check{w}_0(1) X +
     \rO^*\Bigl(27\zeta(1+\varepsilon)^{3n_{\K}} |d_{\K}|^{\frac{1+\varepsilon}{2}}\N\q^{\frac{1+\varepsilon}{4}}
     \theta(\q)^2X^{\frac{1+\varepsilon}{2}}M^*(\varepsilon,1)\Bigr),
 \end{align*}
  where $\alpha_{\K, \q} = \frac{\alpha_{\K}\phi(\q)J(1, \chi_{0, \q})}{\N\q}$ 
  is the residue of $F(s, \chi_{0, \q})$ at $s=1$ and $M^*$ is
 as defined in Lemma~\ref{getMstar}. We also have
 \begin{equation*}
    \alpha_{\K, \q} = \alpha_{\K} J(1, \chi_{0, \q})\frac{\phi(\q)}{\N\q}
    ~\ge~
     \frac{\phi(\q)}{\N\q}\frac{\alpha_{\K}}{\sqrt{\zeta(2)^{n_{\K}}}}
  \end{equation*}
where we have used the inequality~\eqref{lb-j}. The
above is valid for a general smoothing function~$w$ but we restrict ourselves to $w = w_0$. 

Applying Lemma~\ref{studyw0} and using yet again the notation $n$ for $n_\K$,
we get 
\begin{equation*}
  S_1(w_0)
 ~\ge~
  \frac{9(4/\zeta(2))^nX}{500\cdot
    \sqrt{n\N\q}~|d_\K|}\frac{\phi(\q)}{\N\q}
  ~-~
  324~\zeta(1+\varepsilon)^{3n} (57n)^{n + 3} |d_{\K}|^{\frac{1+\varepsilon}{2}}\N\q^{\frac{1+\varepsilon}{4}}
     \theta(\q)^2X^{\frac{1 +\varepsilon}{2}} 
\end{equation*}
and the first summand is larger than twice the second one provided
that
\begin{equation*}
  X^{\frac{1-\varepsilon}{2}}
  ~\ge~ 7 \cdot 10^9 n^{\frac{7}{2}} (24n)^n\zeta(1+\varepsilon)^{3n}
  |d_{\K}|^{\frac{3+\varepsilon}{2}}
  \theta(\q)^{2}
   \frac{\N\q}{\phi(\q)}\N\q^{\frac{3+\varepsilon}{4}} .
 \end{equation*}
 We select $\varepsilon=1/10$. The previous inequality is implied by
  \begin{equation*}
    X ~\ge~  8 \cdot 10^{21}\cdot n^{\frac{70}{9}} (31944 n)^{\frac{20n}{9}}
    |d_\K|^{\frac{31}{9}} \left(\theta(\q)^{2} \frac{\N\q}{\phi(\q)} \right)^{\frac{20}{9}} \N\q^{\frac{31}{18}}.
  \end{equation*}
Since $\theta(\q)^{2} \frac{\N\q}{\phi(\q)} = \prod_{\P | \q} \frac{\N\P^2}{(\N\P -1)(\sqrt{\N\P} -1)^2}$
and 
$$
\left(\frac{p^2}{(p-1)(\sqrt{p}-1)^2} \right)^{20/9} ~\le~ p^{1/6} 
\text{ when } p > 56 \phantom{m}\text{and}\phantom{m}
 \prod_{p\le 56}\frac{1}{p^{1/6}}\left(\frac{p^2}{(p-1)(\sqrt{p}-1)^2} \right)^{20/9}
  ~\le~ 9 \cdot 10^9,
  $$
  we have
  $$
  X ~\ge~  8 \cdot (10^{31} n^7)^{n}  |d_\K|^{\frac{31}{9}}  \N\q^{\frac{34}{18}}.
  $$
 We simplify the final statement by noticing that
$\frac{9(4/\zeta(2))^n}{1000\sqrt{n}}\ge \frac{1}{27}$.
\end{proof}
%%%%%%%%%%%

We may now complete the proof of Theorem~\ref{primeinkernel}.
%%%%%%%%%%%%%
\begin{proof}[Proof of Theorem~\ref{primeinkernel}]
Suppose that the theorem is not true. Then every degree one prime ideal $\P$ 
co-prime to $\q$ with norm at most 
$
X=
\N\P \le 8 \cdot (10^{31} n_\K^7)^{n_\K}  |d_\K|^4 \N\q^2
$
satisfies the property that $\chi(\P)=-1$. Then for every non square-full ideal 
$\mathfrak{a} ~(\neq \mathcal{O}_{\K})$ which decomposes only as a product 
of prime ideals of degree one and of norm at most $X$, we have
$(1\star \chi)(\mathfrak{a})=0$. By Lemma~\ref{simAxer}, we get a
contradiction, and
this completes the proof of the theorem.
\end{proof}

%%%%%%%%%%%%%%%%%%%%%%%%%%%%%%%%%%%%%%%%%%%%%%%%%
%%%%%%%%%%%%%%%%%%%%%%%%%%%%%%%%%%%%%%%%%%%%%%%%%
\section{Selberg sieve for number fields in sieve dimension one}
%%%%%%%%%%%%%%%%%%%%%%%%%%%%%%%%%%%%%%%%%%%%%%%%%%
%%%%%%%%%%%%%%%%%%%%%%%%%%%%%%%%%%%%%%%%%%%%%%%%%%
\label{Ssnfdone}
In this section, we  derive some lemmas which are required
to prove the number field analogues of two versions of the Brun-Titchmarsh Theorem.
Throughout this section, $z$ will denote a real number greater than one and all ideals
considered are integral ideals. For a fixed integral ideal $\q$, we define 
\begin{equation}\label{prod-def}
V(z) 
= 
\prod_{\P \mid \mathcal{P}(z) \atop{(\P, \q)=\O_\K} } \P
\phantom{mmm}\text{where} \phantom{mmm}
\mathcal{P}(z) 
=
\prod _{\P \atop{\N(\P) \le z} } \P .
\end{equation}
Recall the definition of the M\"{o}bius function in~\eqref{defmoebius}.
Further, for any ideal $\e$ of $\O_{\K}$, we define
\begin{equation}
G_{\e}(z) 
=
\sum_{0< \N(\a) \le z, \atop{(\a, \e) = \O_{\K}}} \frac{\mu^2(\a)}{\phi(\a)},
\quad
G(z)
= 
G_{\O_{\K}}(z).
\end{equation}
For some fixed integral ideal $\q$, we further set
\begin{equation}
  \label{deflambdae}
  \lambda_{\e}(\q)
  =
  \mu(\e) \frac{\N(\e)G_{\e\q}(\frac{z}{\N(\e)})}{\phi(\e)G_\q(z)}
  ,\quad
  \lambda_\e=\lambda_\e(\O_\K).
\end{equation}
We set $\lambda_{\e}(\q) = 0$ whenever  $\N(\e) > z$ or $(\e, \q) \ne \O_{\K}$.
In this set-up, we recall the following two lemmas from \cite{Schaal*68} by
Schaal. The reader may also refer to the beginning of Section~4 of
the paper by Debaene \cite{Debaene*19}.
%%%%%%%%%%%
\begin{lem}
For any ideal $\e \mid V(z)$, one has $|\lambda_{\e}(\q)| \le 1$.
\end{lem}
%%%%%%%%%%%
%%%%%%%%%%%
\begin{lem}\label{G(z)bound}
We have
$
\displaystyle G_\q(z)^{-1} 
= 
\sum_{\e_1, \e_2 \atop{\e_i \mid V(z)}} \frac{\lambda_{\e_1}(\q)  \lambda_{\e_2}(\q)}{\N[\e_1,\e_2]}
$.
\end{lem}
%%%%%%%%%%%%

%%%%%%%%%%%%%
\begin{proof}
  In \cite{Schaal*68}, the required definitions are in $(1.25)$ and in
  $(3.1)$ except
  that we do not consider a more severe sieving restriction, so that
  $\rho=x$ and this latter quantity is called~$z$ in our setting.
  The bound for $|\lambda_{\e}(\q)|$ is given in $(3.2)$. The expression
  for $G_\q(z)$ is contained in the last displayed equations at the
  bottom of page~293 with the remark of eq.~$(3.3)$.
\end{proof}
%%%%%%%%%%%%%

%%%%%%%%%%%%%%%%%%%%%%%%%%%%%%%%%%%%%%%%%%%%%
\subsection{A lower bound for $G_\q(z)$}
%%%%%%%%%%%%%%%%%%%%%%%%%%%%%%%%%%%%%%%%%%%%%
%%%%%%%%%%%
\begin{lem}
  \label{Mellin1}
  When $y>0$ and $k$ is a positive integer, we have
  \begin{equation*}
    \max(0, 1-y)^k = \frac{1}{2\pi i}
    \int_{\Re s=2} y^{-s} \frac{k! ds}{s(s+1) \cdots (s+k)}.
  \end{equation*}
\end{lem}
%%%%%%%%%%%

The following theorem gives a lower bound  for $G_\q(z)$.
%%%%%%%%%%%%%%%%%%%%%%%%%%%
\begin{thm}\label{Gzbound}
When $z\ge (10^6n_\K)^{4 n_\K}|d_\K|^3$, we have
$\displaystyle
G_\q(z)
\ge 
\alpha_{\K} \frac{\phi(\q)}{\N(\q)}
\log\frac{z}{e^2 n_\K\sqrt{|d_{\K}|}}$.
\end{thm}
%%%%%%%%%%%%%%%%%%%%%%%%%%%%
This is a version of Lemma~5 of \cite{Schaal*70} where the dependence
in the field is explicit. In Lemma 14 of \cite{Debaene*19}, a similar result
is proved, but it relies on $R_\K h_\K$ when we prefer to rely on~$d_\K$.

\begin{proof}
  We first remove the dependence in $\q$ by using the following
  inequality from \cite{Schaal*70}, page~266, obtained by combining
  the points (a) and (b) therein:
  \begin{equation}
    \label{Halberstam-Schaal}
    G_\q(z)\ge
    \frac{\phi(\q)}{\N(\q)}
    \sum_{0< \N(\a) \le z} \frac{1}{\N(\a)}.
  \end{equation}
  We note in passing that it is
  straightforward to adapt the inequality \cite[(1.3)]{van-Lint-Richert*65}
  by Van Lint and Richert to prove that $G_\q(z)\ge
    \frac{\phi(\q)}{\N(\q)}G(z)$, which is
    more refined than \eqref{Halberstam-Schaal}.
To handle right hand side of \eqref{Halberstam-Schaal}, we  use Lemma~\ref{Mellin1} with 
$y = {\mathfrak{N}(\a)}/{z}$ to get
\begin{align*}
\sum_{0< \N(\a) \le z} \frac{1}{\N(\a)} 
& \ge 
\sum_{\a \atop \a \ne0} \frac{1}{\N(\a)} \max\left( 0,  
~\left(1 - \frac{\mathfrak{N}(\a)}{z} \right)\right)^k \\
& = 
 \frac{1}{2\pi i}  \sum_{\a \atop \a \ne 0} \frac{1}{\N(\a)} 
 \int_{\Re s=2} \left(\frac{\mathfrak{N}(\a)}{z}\right)^{-s} \frac{k! }{s(s+1) \cdots (s+k)}~ds.
\end{align*}
Since we are in the region of absolute convergence of $\zeta_{\K}(s)$,
this leads to
\begin{equation*}
\sum_{0< \N(\a) \le z} \frac{1}{\N(\a)} 
\ge  
\frac{1}{2 \pi i} \int_{\Re s=2} \zeta_{\K}(1+s)  ~\frac{k! z^s  }{s(s+1) \cdots (s+k)} ~ds.
\end{equation*}
We note that the integrand has a double pole at $s=0$ and move the
line of integration to $\Re s = -1/4$. In the neighborhood of $s=0$, we find that
\begin{equation*}
  s^2\zeta_{\K}(1+s)  \frac{k!  z^s}{s(s+1) \cdots (s+k)}
  =\frac{\alpha_\K k!  z^s}{(s+1) \cdots (s+k)}
  +\frac{\alpha_\K\gamma_\K k! s\, z^s}{(s+1) \cdots (s+k)}
  + O(s^2).
\end{equation*}
The residue at $s=0$ is thus given by
\begin{equation*}
r = \lim_{s \to 0} \frac{d}{ds} \left(s^2 \zeta_{\K}(1+s)  ~\frac{k!  z^s}{s(s+1) \cdots (s+k)}\right) 
=
 \alpha_{\K}  \left( \gamma_{\K}  + \log z 
 -\sum_{\ell=1}^{k} \frac{1}{\ell} \right).
\end{equation*}
The remaining integral is now
\begin{equation*}
  I = \frac{1}{2 \pi i} \int_{\Re s = -\frac{1}{4}}
  \frac{k!  \zeta_{\K}(1+s) z^s}{s(s+1)\cdots(s+k)} ds.
\end{equation*}
Further Lemma~\ref{HR1} with $\varepsilon=1/4$ gives us
$$
|\zeta_{\K} (3/4 + it ) |
\le
27 \zeta(5/4)^{n_{\K}} (|d_{\K}|(7/4 + |t|)^{n_{\K}})^{1/4}.
$$
This gives 
\begin{equation*}
| I |
 \le  
5  ~\zeta(5/4)^{n_{\K}} k! 
\left(\frac{|d_{\K}|}{z}\right)^{1/4}
\int_{\Re s  = -\frac{1}{4}} \frac{(7/4 + |t| )^{\frac{n_{\K}}{4} }|ds|}{ |s(s+1)\cdots(s+k)|}.
\end{equation*}
To estimate the integral, say $I_0$, choose $T>0$, $k>n_\K/4$ and write
\begin{align*}
\tfrac12 I_0
& = 
\int_{t=0, \Re s = -\frac{1}{4} }^T  
\frac{ (7/4 + t )^{ \frac{n_{\K}}{4} } |ds|}{ |s(s+1)\cdots(s+k)|}
 ~+~
\int_{t=T,  \Re s = -\frac{1}{4} }^{\infty} 
               \frac{ ( 7/4 + t )^{ \frac{n_{\K}}{4} } |ds|}{ |s(s+1)\cdots(s+k)|}
  \\
& \le  
\frac{ 8 }{(k-1)!} \int_{0}^T 
\left( \frac{7}{4}+ t \right)^{\frac{n_{\K}}{4}} dt 
~+~
\int_{T}^{\infty} \Bigl(\frac7{4t} + 1\Bigr)^{\frac{n_{\K}}{4}}
       \frac{dt}{t^{k+1-\frac{n_\K}{4}}}
  \\&\le
  \frac{ 32 }{(k-1)!(4+n_\K)}(7/4 + T)^{1+\frac{n_{\K}}{4}}
  ~+~
  \Bigl(\frac7{4T} + 1\Bigr)^{\frac{n_{\K}}{4}}
  \frac{1}{(k-\frac{n_\K}{4})T^{k-\frac{n_\K}{4}}}.
\end{align*}
Hence, with $T=k/e$,
\begin{align*}
  k! I_0
  &\le
  (7/4 + k/e)^{1+\frac{n_{\K}}{4}}
  \biggl(
  \frac{ 64 k}{4+n_\K}  ~+~ \frac{8\cdot k!}{(4k-n_\K)(k/e)^k}
  \biggr)
  \\&\le
  (7/4 + k/e)^{1+\frac{n_{\K}}{4}} \biggl(\frac{ 64 k}{4+n_\K}
  ~+~\frac{8\cdot \sqrt{2\pi k}e^{\frac1{12k}}}{(4k-n_\K)} \biggr)
\end{align*}
by using the explicit Stirling Formula
  recalled in \eqref{ExplicitStirling}.
Therefore, with the choice $k  = n_{\K}$, we find that
\begin{align*}
| I | 
& \le  
\biggl(5^{1/n_\K}
   \zeta(5/4)\biggl(\frac{7}{4n_\K}+\frac{1}{e}\biggr)^{\frac{1}{n_\K}+\frac{1}{4}}
   \biggl(\frac{ 64 n_\K}{4+n_\K}
 +\frac{8\cdot \sqrt{2\pi}e^{\frac1{12n_\K}}}{3\sqrt{n_\K}}\biggr)^{1/n_\K}
   n_\K^{\frac{1}{n_\K}+\frac{1}{4}} \biggr)^{n_{\K}} 
   \left(\frac{ |d_{\K}|}{z}\right)^{\frac{1}{4}}
 \\
&\le  
\bigl(10^{3}n_\K\bigr)^{n_\K}\left(\frac{ |d_{\K}|}{z}\right)^{\frac{1}{4}}.
\end{align*}
Whence
\begin{align*}
  \frac{\N\q}{\phi(\q)} G_\q(z)
  &\ge \alpha_{\K}  \left( \gamma_{\K} + \log z 
    -\log n_\K-1 \right)
  -\bigl(10^{3}n_\K\bigr)^{n_\K}\left(\frac{
    |d_{\K}|}{z}\right)^{\frac{1}{4}}
  \\&\ge
  \alpha_{\K}  \left( \log z -\log\sqrt{|d_\K|}
    -\log n_\K-1 
  -\bigl(10^{3}n_\K\bigr)^{n_\K}\frac{100 \sqrt{|d_{\K}|}}{9 \cdot 2^{n_{\K}}}\left(\frac{
    |d_{\K}|}{z}\right)^{\frac{1}{4}}\right)
\end{align*}
by Lemma~\ref{boundalphaK} and the inequality $\gamma_\K\ge
-\log\sqrt{|d_\K|}$ from \eqref{IharaIneq}.
This completes the proof of \thmref{Gzbound}.
\end{proof}

%%%%%%%%%%%%%%%%%%%%%%%%%%%%%%%%%%%%%%%%%%
\subsection{Controlling the error term in Selberg's sieve}
%%%%%%%%%%%%%%%%%%%%%%%%%%%%%%%%%%%%%%%%%%

Our first lemma borrows from~\cite{Halberstam-Richert*79} by
Halberstam and Richert.

%%%%%%%%%%%%%%%%%%%
\begin{lem}\label{primes}
When $x\ge1$ and $\alpha\in[0,1)$, we have
$\displaystyle
  \sum_{\N\P\le x}\frac{\log\N\P}{\N\P^\alpha}\le
  \frac{1.02 \,n_\K x^{1-\alpha}}{1-\alpha}$.
\end{lem}
%%%%%%%%%%%%%%%%%%%

%%%%%%%%%%
\begin{proof}
Let $\P_1,\cdots, \P_g$ be prime ideals of $\O_\K$ lying over
a rational prime $p$.  We then have
$$
\displaystyle\sum_{1\le i\le g}e_i\log \N\P_i
=
n_\K\log p,
$$
where $e_i\ge1$ is the ramification index of $\P_i$ above~$p$. This
  implies that $\sum_{1\le i\le g}\log \N\P_i\le n_\K\log
  p$. Inequality~\eqref{rosserpsi}
  concludes the proof of the lemma
  when $\alpha=0$. For the general case, we use summation by parts to write
  \begin{align*}
    \sum_{\N\P\le x}\frac{\log\N\P}{\N\P^\alpha}
    &=
      \alpha\int_1^x \sum_{\N\P\le t}\log\N\P\frac{dt}{t^{1+\alpha}}
      +\sum_{\N\P\le x}\frac{\log\N\P}{x^\alpha}
    \\&\le~
      1.02
    \,n_\K\left(\alpha\int_1^x\frac{dt}{t^{\alpha}}+x^{1-\alpha}\right)
    ~\le~ 1.02
    \,n_\K\frac{x^{1-\alpha}}{1-\alpha}
  \end{align*}
  as required.
\end{proof}
%%%%%%%%%%
%%%%%%%%%%%%%%%%%%%
\begin{lem}\label{userhhr}
When $x\ge 1$ and $\alpha\in[0,1)$, we have
\begin{equation*}
  \sum_{0< \N\a\le x}\mu^2(\a)\left(\frac{x}{\N\a}\right)^\alpha
  ~\le~
  \frac{(1+1.02 n_\K)x}{(1-\alpha)(1+\log x)}\sum_{0< \N\a\le x}\frac{\mu^2(\a)}{\N\a}.
\end{equation*}
\end{lem}
%%%%%%%%%%%%%%%%%%%

%%%%%%%%%%%%%%
\begin{proof}
  We set $S(x)=\sum_{ 0<\N\a\le x}\mu^2(\a)\left({x}/{\N\a}\right)^\alpha$.
  We use $\log y\le (y^{1-\alpha}-1)/(1-\alpha)$ when $y>0$ and
  readily find that
  \begin{align*}
    S(x)\log x
    &=\sum_{0< \N\a\le x}\mu^2(\a)\left(\frac{x}{\N\a}\right)^\alpha\log\frac{x}{\N\a}
       +\sum_{0< \N\a\le x}\mu^2(\a)\left(\frac{x}{\N\a}\right)^\alpha\log\N\a
    \\&\le
    \frac{x}{1-\alpha}\sum_{0< \N\a\le x}\frac{\mu^2(\a)}{\N\a}-S(x)
    +\sum_{0< \N\a\le x}\mu^2(\a)\left(\frac{x}{\N\a}\right)^\alpha\sum_{\P|\a}\log\N\P
    \\&\le
    \frac{x}{1-\alpha}\sum_{0< \N\a\le x}\frac{\mu^2(\a)}{\N\a}-S(x)
    +\sum_{\N\P\le x}\log\N\P
    \sum_{0< \N\a\le x/\N\P} \mu^2(\a)\left(\frac{x/\N\P}{\N\a}\right)^\alpha.
  \end{align*}
  We invert the summation in the last term and
  then appeal to Lemma~\ref{primes} to get
  \begin{equation*}
  S(x)(1+\log x)
  ~\le~
    \frac{x}{1-\alpha}\sum_{0< \N\a\le x}\frac{\mu^2(\a)}{\N\a}
    +\frac{1.02 n_\K\,x }{1- \alpha} \sum_{0< \N\a\le x}\frac{\mu^2(\a)}{\N\a}
  \end{equation*}
  and our lemma follows readily.
\end{proof}
%%%%%%%%%%%%%%

%%%%%%%%%%%%%%%%%%%%%%
\begin{thm}\label{absolute}
When $z\ge1$ and $\alpha\in[0,1)$, we have
  \begin{equation*}
    \sum_{\substack{\e \mid V(z)}} \frac{|\lambda_{\e}(\q)|}{\N\e^\alpha}
    \le
    \frac{3.1\,n_\K z^{1-\alpha}}{(1-\alpha)(2+\log z)}
    c_1(\alpha)^{n_\K}
    +
    z^{3(1-\alpha)/4}c_2(\alpha)^{n_\K},
  \end{equation*}
  where
  \begin{equation*}
  c_1(\alpha)=
  \prod_{p}\biggl(1+\frac{1+p^\alpha}{(p-1)p}\biggr),
  \quad
  c_2(\alpha)=
   \prod_{p}\biggl(1+\frac{1+p^\alpha}{(p-1)p^{\frac{1+ 3\alpha}4}}\biggr).
\end{equation*}  
\end{thm}
%%%%%%%%%%%%%%%%%%%%%%
This theorem replaces Lemma~15 in \cite{Debaene*19}, but notice that we avoid the
high power of $\log z$ and even save an additional $\log z$. We also
mention for future use that
\begin{align}\label{boundc2}
  \frac{c_2(\alpha)}{\zeta(\frac{5+3\alpha}{4})\zeta(\frac{5-\alpha}4)}
  ~=~
  \prod_{p}\biggl( 1 +  \frac{(p^{1/4} -p^{\alpha/4})p^{\alpha/4}}{p^{5/2} - p^{3/2}}
    ~+~ \frac{ -(p^{5/4} + p^{1/4}) p^{\alpha} + p^{\frac{6 + 3\alpha}{4} }  + p^{ \frac{\alpha}{4} } 
    + p^{\frac{5\alpha} {4} } - p^{5/4}}{ p^{\frac{3\alpha}{  2} } (p^{15/4} - p^{11/4})}   
  \biggr).
\end{align}
where the right hand side is uniformly bounded above for $\alpha\in[0,1]$.
%%%%%%%%%%%%%
\begin{proof}
Set $T = \sum_{\e \mid V(z)} |\lambda_{\e}(\q)|/\N\e^\alpha$. Then from the definition of $\lambda_{\e}(\q)$, we get
\begin{align*}
  G_\q(z)T
  &=  
    \sum_{\e \mid V(z)} \frac{\N\e^{1-\alpha}}{\phi(\e)}  
    \sum_{\N\a \le \frac{z}{\N\e} ,\atop{(\a, \e\q) = \O_{\K}}} \frac{\mu^2(\a)}{\phi(\a)}
    = 
    \sum_{\e \mid V(z)} \N\e^{1-\alpha}
    \sum_{\N\b \le z, \atop{(\b,\q)= \O_{\K}, \atop{\e \mid \b}}} 
    \frac{\mu^2(\b)}{\phi(\b)}
  \\&=
  \sum_{\N\b \le z, \atop{(\b,\q)= \O_{\K}}}
  \frac{\mu^2(\b)}{\phi(\b)} \sum_{\e \mid \b} \N\e^{1-\alpha} .
\end{align*}
We note that, when $\b$ is squarefree, we have
$$
 \frac{\N\b^\alpha}{\phi(\b)} \sum_{\e \mid \b} \N\e^{1-\alpha}
 = 
  \prod_{\P \mid \b}
 \left(
   \frac{\N\P + \N\P^{\alpha}}{\N\P-1}\right) 
 =
  \prod_{\P \mid \b} \left(1+\frac{1+\N\P^\alpha}{\N\P-1}\right) 
=
 \sum_{\f\mid \b}  g_\alpha(\f)
$$
where $g_\alpha(\f)$ is  the completely multiplicative function defined on primes
by $g_\alpha(\P)=\frac{1+\N\P^\alpha}{\N\P-1}$. Only the values on
squarefree ideals are required but to define a unique function, we
added that it is a `completely' multiplicative function. 
Therefore 
\begin{align*}
  G_\q(z)T 
  &= 
    \sum_{\N\b \le z, \atop{(\b,\q) = \O_{\K}}}
    \frac{\mu^{2}(\b)}{\N\b^\alpha}
    \sum_{\f \mid \b}  g_\alpha(\f)
    = 
    \sum_{\N\f \le z} \mu^2(\f) g_\alpha(\f)
    \sum_{\N\b \le z, \atop{(\b,\q)=\O_{\K}, \atop{\f \mid \b}}}
    \frac{\mu^{2}(\b)}{\N\b^\alpha} \\
  &\le
    \sum_{\N\f \le z, \atop{(\f,\q) = \O_{\K}}}
    \frac{\mu^2(\f) g_\alpha(\f)}{\N\f^\alpha}
    \sum_{\N\b_1 \le \frac{z}{\N\f}}\frac{\mu^{2}(\b_1)}{\N\b_1^\alpha}.
\end{align*}
At this level, we use the trivial bound 
$$
G_\q(z) = \sum_{\N(\a) \le z \atop{(\a, \q)= \O_{\K}}} 
\frac{\mu^2(\a)}{\N\a} \prod_{\P | \a} \sum_{k \ge 0} \N\P^{-k}
~\ge  
\sum_{\N(\a) \le z, \atop{(\a, \q)= \O_{\K}}} \frac{1}{\N(\a)}
$$
and Lemma~\ref{userhhr} with $x=z/\N\f$ for the inner sum
of the right-hand side when $\N\f\le \sqrt{z}$, and with $x=\sqrt{z}$
for the remaining sum. This first gives us

\begin{multline*}
  \sum_{\N(\a) \le z, \atop{(\a, \q)= \O_{\K}}} \frac{1}{\N(\a)} T
  ~\le~
  \frac{(1+1.02n_\K)z^{1-\alpha}}{1-\alpha}
  \sum_{\N\f \le \sqrt{z}, \atop{(\f,\q) = \O_{\K}}}
  \frac{\mu^2(\f) g_\alpha(\f)}{\N\f(1+\log(z/\N\f))} 
  \sum_{\N\b_1 \le \frac{z}{\N\f}, \atop{ (\b_1, \q)= \O_{\K}}}\frac{\mu^2(\b_1)}{\N\b_1}
  \\+~
  z^{(1-\alpha)/2}\sum_{z \ge \N\f > \sqrt{z}, \atop{(\f,\q) = \O_{\K}}}
  \frac{\mu^2(\f) g_\alpha(\f)}{\N\f^\alpha}
  \sum_{\N(\a) \le \sqrt{z}, \atop{(\a, \q)= \O_{\K}}} \frac{1}{\N(\a)}.
\end{multline*}
We readily deduce from this the inequality
\begin{equation*}
   T
  \le
  \frac{3.1\,n_\K z^{1-\alpha}}{(1-\alpha)(2+\log z)}
  \sum_{\N\f \le z, \atop{(\f,\q) = \O_{\K}}}
  \frac{\mu^2(\f)g_\alpha(\f)}{\N\f}
  ~+~ 
  z^{(1-\alpha)/2}\sum_{z \ge \N\f > \sqrt{z}, \atop{(\f,\q) = \O_{\K}}}
   \frac{\mu^2(\f) g_\alpha(\f)}{\N\f^\alpha}.
\end{equation*}
To proceed, we forget
about the coprimality in~$\f$ and use Rankin's trick for the second
part, noticing that $(z/\N\f)^{(1-\alpha)/4}\ge 1$ therein. This leads~to
\begin{equation*}
  T
  \le
  \frac{3.1\,n_\K z^{1-\alpha}}{(1-\alpha)(2+\log z)}C_1(\alpha)
  ~+~ z^{3(1-\alpha)/4}C_2(\alpha).
\end{equation*}
Concerning the constants $C_1(\alpha)$ and $C_2(\alpha)$, we readily find that
\begin{equation*}
  \left\{
    \begin{aligned}
  C_1(\alpha)&=\sum_{\f \atop \f \ne 0}\frac{\mu^2(\f) g_\alpha(\f)}{\N\f}
 ~ =~
 \prod_{\P}\biggl(1+\frac{1+\N\P^\alpha}{(\N\P-1)\N\P}\biggr)
  ~\le~
  \prod_{p}\biggl(1+\frac{1+p^\alpha}{(p-1)p}\biggr)^{n_\K},
  \\
  C_2(\alpha)
  &=~
  \sum_{\f \atop \f \ne 0}\frac{\mu^2(\f) g_\alpha(\f)}{\N\f^{\frac{3 \alpha +1}{4}}}
  ~=~
  \prod_{\P}\biggl(1+\frac{1+\N\P^\alpha}{(\N\P-1)\N\P^{\frac{3\alpha+1}4}}\biggr)
  ~\le~ 
  \prod_{p}\biggl(1+\frac{1+p^\alpha}{(p-1)p^{\frac{3\alpha+1}4}}\biggr)^{n_\K}.
\end{aligned}
\right.
\end{equation*}
The lemma follows readily.
\end{proof}
%%%%%%%%%%%%%%%
%%%%%%%%%%%
\begin{cor}
  \label{absolute0}
  When $z\ge1$, we have
  $\displaystyle
    \sum_{\substack{\e \mid V(z),\\ (\e,\q)=\O_\K}} |\lambda_{\e}(\q)|
    \le
    6\cdot
    89^{n_\K}\frac{z}{2+\log z}$.
\end{cor}
%%%%%%%%%%%

%%%%%%%%%%
\begin{proof}
We use Theorem~\ref{absolute} with $\alpha=0$. We numerically find that
\begin{equation*}
  1.6\,n_\K \prod_{p}\biggl(1+\frac{2}{(p-1)p}\biggr)^{n_\K}
  ~\le~
  \prod_{p}\biggl(1+\frac{2}{(p-1)p^{1/4}}\biggr)^{n_\K}\le (88.2)^{n_\K}.
\end{equation*}
Furthermore $z^{3/4}\le 4\,z/(2+\log z)$ and $6\,(88.2)^{n_\K}\le 6\cdot
89^{n_\K}$.
\end{proof}
%%%%%%%%%%
%%%%%%%%%
\begin{cor}
  \label{absolute1minus1overnK}
  When $z\ge1$, we have
  $\displaystyle
    \sum_{\substack{\e \mid V(z),\\ (\e,\q)=\O_\K}}
    \frac{|\lambda_{\e}(\q)|}{\N\e^{1-\frac{1}{n_\K}}}
    ~\le~
    n_\K^{9n_\K}\frac{z^{1/n_\K}}{2+\log z}$.
\end{cor}
%%%%%%%%%%%%%%

%%%%%%%%%%%%%%
\begin{proof}
  We use Theorem~\ref{absolute} with $\alpha=1-1/n_\K$. Let us call $T$ the quantity to be bounded
  above. Theorem~\ref{absolute} gives the upper bound, with~$n=n_\K$,
  \begin{multline*}
    \frac{3.1\,n^2 z^{1/n}}{2+\log z}
    c_1(1-1/n)^{n}
    +
    z^{3/(4n)}c_2(1-1/n)^{n}
    \\=
    \frac{n\, z^{1/n}}{2+\log z}
    \Bigl(3.1\,n c_1(1-1/n)^{n} + 4c_2(1-1/n)^{n}\frac{\frac1{2n}+\log(z^{1/(4n)})}{z^{1/(4n)}}\Bigr).
  \end{multline*}
  The maximum of $y\mapsto(\frac14+\log y)/y$ is attained at $y=e^{3/4}$ and
  $c_1(\alpha)\le c_2(\alpha)$, so we can simplify this upper bound to
  \begin{equation*}
    \frac{n\, z^{1/n}}{2+\log z}(3.1\,n +1.9)c_2(1-1/n)^{n}.
  \end{equation*}
  Now we want to find an upper bound for $c_2(1-1/n)$. 
  Applying inequality~\eqref{boundc2} and the fact that $\alpha\ge1/2$
  in our case, we get
  \begin{equation*}
    \frac{c_2(1-1/n)}{\zeta(1+ \frac{1}{4n})\zeta(2-\frac{3}{4n})}
    ~\le~ 
    \prod_{p}\biggl( 1 +  \frac{p^{1/2}}{p^{5/2} - p^{3/2}}
    + \frac{p^2 - p^{3/2}  - p^{1/2} + 1}{p^{1/2}(p^{15/4} - p^{11/4})}
    \biggr).
  \end{equation*}
A direct computation shows that right hand side of the product for $p\le 10^4$ is bounded by
$2.3$, and then using calculus 
 we derive that ${c_2(1-1/n)}\le 11.5n \zeta(13/8)$. This leads
  to
  \begin{equation*}
    c_2(1-1/n) ~\le~  26.45\,n.
  \end{equation*}
  Furthermore the quantity $n\, (3.1\,n +1.9)(26.45\,n)^{n}$ may further be bounded above by
  $n^{9n}$. This establishes this lemma.
\end{proof}
%%%%%%%%%%%%%%

%%%%%%%%%%%%%%%%%%%%%%%%%%%%%%%%%%%%%%%%%%%%%%%%%%%%%%%%%%%%%%%%%%%%%%%%%%%%%%%%%%%%%%%%%
\section{Brun-Titchmarsh Theorem for cosets of ray class groups}
%%%%%%%%%%%%%%%%%%%%%%%%%%%%%%%%%%%%%%%%%%%%%%%%%%%%%%%%%%%%%%%%%%%%%%%%%%%%%%%%%%%%%%%%%
\label{btc}
In this section, we prove a number field analogue of Brun-Titchmarsh theorem for
cosets. This result, while being of independent interest, is also crucial in proving our 
main theorem. 
%%%%%%%%%%%%%%%%%%%%%%
\begin{thm}\label{bt}
Let $\a, \q$ be integral ideals with $(\a, \q) = \O_{\K}$. Also let
$H$ be a subgroup of $H_{\q}(\K)$ with index $Y$. When 
$X/Y
~\ge~
(10^6n_\K)^{8n_\K + 11}|d_\K|^6\sqrt{|d_{\K}|\N(\q)}\log(|d_\K|\N\q)^{n_\K}$, 
we have
\begin{equation*}
  \sum_{\P \in [\a]H}  w\left(\frac{\N(\P)}{X}\right) 
  \le 
  \frac{2 \|w\|_1 X}
  {Y\log\frac{u^*(w,\K)X}{Y\sqrt{\N(\q)}\log(|d_\K|\N\q)^{n_\K}}}~,
  \quad%
  u^*(w,\K)=\frac{\|w\|_1/(\|w^{(n_\K+3)}\|_\infty+5\|w\|_1)}{20000\cdot2^{22n_\K}|d_{\K}|^{3/2}}.
\end{equation*}
  for any \emph{non-negative} smoothing function $w$ as defined in
  Section~\ref{Notation}.
\end{thm}
%%%%%%%%%%%%%%%%%%%%%
When compared to the usual Brun-Titchmarsh Theorem, this result has
the surprising feature to be sometimes valid for instance when
$X<\N\q$. For instance, when $Y=2$, the bound $X$ can be as small as
$\rO_{\K}(\N\q^{\frac12+\epsilon})$. 

%%%%%%%%%%%
\begin{proof}
 We use the notation $G=H_{\q}(\K)$ and $n=n_\K$.
For any ideal $\b$ of $\O_\K$, we define
\begin{eqnarray*}
\delta(\b) 
= \frac{1}{Y} \sum_{\chi \in {G\widehat /H }} \chi([\b \a^{-1}]H) 
& \text{ and } & 
\delta^*(\b) 
= \frac{1}{Y} \sum_{\chi \in {G\widehat /H}} \chi^*({[\b \a^{-1}] H})
\end{eqnarray*}
where $G\widehat/H$ denotes the group of characters of $G/H$, 
$\chi^*$ is the primitive character inducing $\chi$. For
any integral ideal $\b$ of $\O_{\K}$ with $(\b, \q) = \O_{\K}$, 
$\overline{\b}$ denotes an element of $G/H$.
In this case, $\delta(\b) = \delta^*(\b)$. For the case $(\b, \q) \neq \O_{\K}$, 
let $\q(\b)$ be the largest divisor of $\q$ co-prime to 
$\b$ and let  $L_{\q(\b)}$ be the image of $H$ in $H_{\q(\b)}(\K) = G_{\q(\b)}$, say.
Note that
$$
\delta^*(\b) 
=
\frac{1}{Y} \sum_{\chi \in {{G_{\q(\b)}\widehat/ L_{\q(\b)}}}} \chi^*({ [\b \a^{-1}]H })
=
\frac{|G_{\q(\b)}/L_{\q(\b)}|}{Y} ~~{\bf 1}_{[\a]L_{\q(\b)}},
$$
where $ {\bf 1}_{ [\a]L_{\q(\b)} }$ is the characteristic function of $[\a]L_{\q(\b)}$.
This shows that $\delta^*(\b)$ is non-negative whenever $(\b, \q) \neq \O_\K$.
Therefore $\delta^*(\b) \ge \delta(\b)$. 

\noindent
Consider the sum $T_1=\sum_{\P} \delta(\P) w({\N(\P)}/{X})$. As
$|w|\le 1$, we readily find that
\begin{align*}
T_1
&= 
      \sum_{\N(\P) \le z} \delta(\P) w\left(\frac{\N(\P)}{X}\right)
      + \sum_{\N(\P) > z}  \delta(\P) w\left(\frac{\N(\P)}{X}\right)
  \\
&\le~
n z + \sum_{(\b, V(z)) = \O_{\K}} \delta^*(\b) w\left(\frac{\N(\b)}{X}\right)
~\le~  
nz  + \sum_{\b} \delta^*(\b) w\left(\frac{\N(\b)}{X}\right) \sum_{\e \mid (\b,V(z))} \mu(\e),
\end{align*}
where $V(z)$ is as in Section~\ref{Ssnfdone}.
Since $\sum_{\e \mid \b}\mu (\e)\le(\sum_{\e \mid \b}\lambda_{\e})^2$, we have
$$
T_1 - nz
~\le~
\sum_{\b} \delta^*(\b) w\left(\frac{\N(\b)}{X}\right) \left(\sum_{\e \mid (\b, V(z))} \lambda_{\e}\right)^2
~\le~ 
\sum_{\e_1, \e_2 \atop{\e_i \mid V(z)}} \lambda_{\e_1} \lambda_{\e_2} 
\sum_{[\e_1,\e_2] \mid \b} \delta^*(\b) w\left(\frac{\N(\b)}{X}\right).
$$
Replacing the definition of $\delta^{*}(\b)$, we get
\begin{align*}
T_1 - nz
& \le~ 
\sum_{\e_1, \e_2 \atop{\e_i \mid V(z)}} \frac{\lambda_{\e_1} 
\lambda_{\e_2}}{Y} \sum_{[\e_1,\e_2] \mid \b} w\left(\frac{\N(\b)}{X}\right) 
                \sum_{\chi \in {G\widehat/H}} \chi^{*}({[\b \a^{-1}]H})
  \\
&\le~ 
\sum_{\e_1, \e_2 \atop{\e_i \mid V(z)}} \frac{\lambda_{\e_1}  \lambda_{\e_2}}{Y} 
\sum_{\chi \in  {G\widehat/H} }\chi^{*}({[\a^{-1}]H}) 
\sum_{[\e_1,\e_2] \mid \b} w\left(\frac{\N(\b)}{X}\right) \chi^*({[\b}]H) 
\\& \le~ 
\sum_{\e_1, \e_2 \atop{\e_i \mid V(z)}} \frac{\lambda_{\e_1}  \lambda_{\e_2}}{Y} 
\sum_{\chi \in {G\widehat/H}}\chi^{*}({[\a^{-1}]H}) \sum_{\m} w\left(\frac{\N(\m[\e_1,\e_2])}{X}\right) 
\chi^*({[\m[\e_1,\e_2]~]H}).
\end{align*}
Applying Mellin transforms, we get 
$$
\sum_{\m} w\left(\frac{\N(\m[\e_1,\e_2])}{X}\right) \chi^*({[\m[\e_1,\e_2]~]H})
 =  
 \chi^{*}([~[\e_1,\e_2]~]H)\frac{1}{2\pi i}\int_{2 - i \infty}^{2+ i\infty} \check{w}(s) L_{\q^*}(s,\chi^*) 
 \frac{X^s ds}{\N([\e_1,\e_2])^s}.
 $$
 where $\q^*$ is the conductor of $\chi^*$. We shift the line of integration to the line $\Re s=0$.
This leads to the upper bound of $T_1 - nz$ given by
$$
\sum_{\e_1, \e_2 \atop{\e_i \mid V(z)}} \frac{\lambda_{\e_1}  \lambda_{\e_2}}{Y} 
\left\{ \frac{\alpha_{\K}\phi(\q)\check{w}(1)X}{\N(\q)\N([\e_1,\e_2])}  
+
 \sum_{\chi \in {G \widehat/H}} \frac{\chi^{*}({[\a^{-1}[\e_1,\e_2]]H})}{2 \pi i}
\int_{-i\infty}^{i\infty} \check{w}(s) L_{\q^*}(s,\chi^*) \frac{X^s ds}{\N([\e_1,\e_2])^s} \right\}.
$$
On using \lemref{HR2} when $\chi^*\neq1$ and \lemref{HR1} when
$\chi^*=1$ on the line $\Re s=0$,  we see that $T_1 - nz$  is less than
$$
\frac{\alpha_\K\check{w}(1) \phi(\q)X}{\N(\q) Y}
\sum_{\e_1, \e_2  \atop{\e_i \mid V(z)}} \frac{\lambda_{\e_1}
\lambda_{\e_2}}{\N([\e_1,\e_2])} 
+
\rO^* \biggl( 3
  M(w,\varepsilon)\zeta(1+\varepsilon)^{n_{\K}}  (|d_{\K}|\N(\q))^{\frac{1+\varepsilon}{2}}
\Bigl(\sum_{\e \atop{\e \mid V(z)}}|\lambda_{\e}|\Bigr)^2 \biggr),
$$
where $0< \varepsilon <1/24$ and
$M(w,\varepsilon)$ is defined and bounded in Lemma~\ref{getM}.
The sum over $\e$ may  be treated by Corollary~\ref{absolute0}.
As a partial conclusion, we reach
\begin{equation*}
  T_1\le n z+ \frac{\alpha_\K\check{w}(1)X}{YG(z)}
  + 12\cdot 2^{\frac{\varepsilon}{2}n}(\|w^{(n+3)}\|_\infty
  +10\cdot 2^{\frac{n}{2}}\|w\|_1)\zeta(1+\varepsilon)^{n}  (|d_{\K}|\N(\q))^{\frac{1+\varepsilon}{2}}36\cdot
89^{2n}\frac{z^2}{\log^2 z}.
\end{equation*}
Let us assume that $z\ge (10^6n)^{4n}|d_\K|^3$ to use
Theorem~\ref{Gzbound}. This leads to
\begin{equation*}
  T_1\le n z+ \frac{\check{w}(1)X}{Y\log\frac{z}{n_\K^4\sqrt{|d_{\K}|}}}
  + 432 \cdot 2^{\frac{28+\varepsilon}{2}n}(\|w^{(n+3)}\|_\infty
  +10\cdot 2^{\frac{n}{2}}\|w\|_1)\zeta(1+\varepsilon)^{n}  (|d_{\K}|\N(\q))^{\frac{1+\varepsilon}{2}}
  \frac{z^2}{\log^2 z}.
\end{equation*}
We select $\varepsilon=\min(1/2, 2/\log(|d_\K|\N\q))$ and use
\begin{equation*}
  \zeta(1+\varepsilon)^{n}  (|d_{\K}|\N(\q))^{\frac{1+\varepsilon}{2}}
  \le e \bigl(2\log(|d_\K|\N\q)\bigr)^{n}\sqrt{|d_{\K}|\N(\q)}
\end{equation*}
to derive the bound
\begin{equation*}
  \frac{Y\log\frac{z}{n^4\sqrt{|d_{\K}|}}}{X \|w\|_1}  T_1 ~ \le
 1
  + 2600\cdot 2^{16n}\cdot \frac{(\|w^{(n+3)}\|_\infty + 5 \|w\|_1)}{  \|w\|_1}
  \log(|d_\K|\N\q)^{n}\sqrt{|d_{\K}|\N(\q)}
  \frac{Y z^2}{X\log z}.
\end{equation*}
We select
\begin{equation*}
  z^2=\frac{X\|w\|_1/Y}{2600\cdot 2^{16n}(\|w^{(n+3)}\|_\infty + 5 \|w\|_1)
  \log(|d_\K|\N\q)^{n}\sqrt{|d_{\K}|\N(\q)}}.
\end{equation*}
This gives us the inequality
\begin{equation*}
  \frac{Y\log\frac{z}{n^4\sqrt{|d_{\K}|}}}{X \|w\|_1}  T_1 
  ~ \le~ 1 +   \frac{1}{\log \frac{z}{n^4 \sqrt{|d_{\K} |} } }
   ~ \le~ \frac{1}{ 1 -   \log^{-1} \frac{z}{n^4 \sqrt{ |d_{\K}| }} },
\end{equation*}
where we have used the inequality $1+x\le 1/(1-x)$ for $x\in[0,1)$.
Thus
\begin{equation*}
   \sum_{\P \in [\a]H}  w\left(\frac{\N(\P)}{X}\right) 
   ~\le~ 
   \frac{2 \|w\|_1 X}{Y\log  \frac{X/ (Y \sqrt{\N\q} \log(|d_\K|\N\q)^{n} ) }{
   |d_{\K}|^{3/2} 20000\cdot 2^{22n}\cdot (\|w^{(n+3)}\|_\infty + 5 \|w\|_1)  / \|w\|_1 }},
 \end{equation*}
 where we have used the inequality $2^{16n}n^8\le 2^{22n}$. This completes the 
 proof of our theorem.
\end{proof}
%%%%%%%%%%%

%%%%%%%%%%%%%%%%%%%%%%%%%%%%%%%%%%%%%%%%%%%%%%%%%%%%%%%%%%%%%%%%%%%%%%%%%%
%%%%%%%%%%%%%%%%%%%%%%%%%%%%%%%%%%%%%%%%%%%%%%%%%%%%%%%%%%%%%%%%%%%%%%%%%%
\section{Brun-Titchmarsh Theorem for single class of ray class groups}
%%%%%%%%%%%%%%%%%%%%%%%%%%%%%%%%%%%%%%%%%%%%%%%%%%%%%%%%%%%%%%%%%%%%%%%%%%
%%%%%%%%%%%%%%%%%%%%%%%%%%%%%%%%%%%%%%%%%%%%%%%%%%%%%%%%%%%%%%%%%%%%%%%%%%
Applying Theorem~\ref{bt} for the trivial subgroup and forgetting
the field-dependence leads to the upper bound
$\frac{2X}{h_{\K,\q}\log(X/\N\q^{3/2+\epsilon})}$ while the usual
Brun-Titchmarsh inequality has $\log(X/\N\q^{1+\epsilon})$. 

%%%%%%%%%%
\begin{proof}[Proof of Theorem~\ref{bt-tri}]
 As before, we denote $n_\K$ by $n$.
We use a Selberg sieve of dimension one as presented in
  Section~\ref{Ssnfdone} with coefficients
  $(\lambda_{\e} (\q))_{\N\e\le z}$.
  On denoting by $T_1$ the sum to evaluate, this gives us
  \begin{equation*}
    T_1
    ~\le~ \sum_{\substack{\a\in [\b],\\ \N\a \le X}}
    \biggl(\sum_{\e \mid (\a, V(z))} \lambda_{\e}(\q)\biggr)^2
    ~+~ n z
    ~\le~
    \sum_{\substack{\e_1,\e_2}}\lambda_{\e_1}(\q)\lambda_{\e_2}(\q)
    \sum_{\substack{[\e_1,\e_2]|\a\in[\b],\\ \N\a\le X}}1 ~+~ nz.
  \end{equation*}
  We write $\a=[\e_1,\e_2]\c$ where $\c$ is an integral ideal
  in the class of $\b[\e_1,\e_2]^{-1}$ (this is legal since the lcm
  $[e_1,\e_2]$ is indeed prime to $\q$), and of norm $\le
  X/\N[\e_1,\e_2]$. We now apply Theorem~\ref{asymfinal}. This gives us
  \begin{multline*}
    T_1-n z
    ~\le~
    \frac{\alpha_{\K} \phi(\q)}{h_{\K,\q}}
    \frac{X}{  \N\q}
    \sum_{\substack{\e_1,\e_2}}
    \frac{\lambda_{\e_1}(\q)\lambda_{\e_2}(\q)}{\N[\e_1,\e_2]}
    +
    \biggl(\sum_{\substack{\e}}
    |\lambda_{\e}(\q)|\biggr)^2
    n^{8n}R_\K F(\q)
    \\+~
    E(\K)
    F(\q)^{\frac{1}{n}}\log(3F(\q))^{n}
    \left( \frac{X}{\N\q}   \right)^{1-\frac{1}{n} }
    \sum_{\substack{\e_1,\e_2}}
    \frac{|\lambda_{\e_1}(\q)\lambda_{\e_2}(\q)|}{\N[\e_1,\e_2]^{1-\frac{1}{n}}}.
  \end{multline*}
  The first term on the right hand side is handled by Lemma~\ref{G(z)bound}, the second
  term is controlled by Corollary~\ref{absolute0} and we now show that the third
  one can be controlled by Corollary~\ref{absolute1minus1overnK}. Indeed,
  Selberg's diagonalisation process gives us, with $\alpha=1-\frac{1}{n}$,
  \begin{equation*}
    \sum_{\substack{\e_1,\e_2}}
    \frac{|\lambda_{\e_1}(\q)\lambda_{\e_2}(\q)|}{\N[\e_1,\e_2]^{\alpha}}
    ~=~
      \sum_{\substack{\e_1,\e_2}}\N(\e_1,\e_2)^{\alpha}
      \frac{|\lambda_{\e_1}(\q)\lambda_{\e_2}(\q)|}{\N\e_1^\alpha\N\e_2^{\alpha}}
   ~ \le~ 
    \sum_{\d, \atop \N\d\le z}\phi_\alpha(\d)
    \biggl(\sum_{\d|\e}
    \frac{|\lambda_{\e}(\q)|}{\N\e^\alpha}\biggr)^2,
  \end{equation*}
  where $\phi_\alpha(\d) = \N(\d)^{\alpha}$.
  We now note that
  \begin{align*}
    \sum_{\d|\e}
    \frac{|\lambda_{\e}(\q)|}{\N\e^\alpha}
    &\le~
      \sum_{(\e,\d\q)=1}
    \frac{\N\e\N\d}{\phi(\e)\phi(\d)\N\d^\alpha\N\e^\alpha}
      \frac{G_{\e\d\q}(\frac{z/\N\d}{\N\e})}{G_\q(z)}
    \\&\le~
    \frac{G_{\d\q}(z/\N\d)}{G_\q(z)}
    \frac{\N\d^{1-\alpha}}{\phi(\d)}
    \sum_{(\e,\d\q)=1}
    \frac{\N\e}{\phi(\e)\N\e^\alpha}
    \frac{G_{\e\d\q}(\frac{z/\N\d}{\N\e})}{G_{\d\q}(z/\N\d)}
    \\&\le~ \frac{G_{\d\q}(z/\N\d)}{G_\q(z)}
    \frac{\N\d^{1-\alpha}}{\phi(\d)}
    \sum_{(\e, \d\q)=1}
    \frac{|\lambda_\e(\q,z/\N\d)|}{\N\e^\alpha},
  \end{align*}
  where $\lambda_\e(\q,z/\N\d) =  \frac{\N\e}{\phi(\e)}
    \frac{G_{\e\d\q}(\frac{z/\N\d}{\N\e})}{G_{\d\q}(z/\N\d)}$. By Corollary~\ref{absolute1minus1overnK},
  this is bounded above by
  \begin{equation*}
    n^{9n}\frac{\N\d}{\phi(\d)}
    \frac{G_{\d\q}(z/\N\d)}{G_\q(z)}\frac{z^{1-\alpha}}{(2+\log(z/\N\d))\N\d}.
  \end{equation*}
 We use the upper bound $G_\q(z)$ for $(\N\d/\phi(\d))G_{\d\q}(z/\N\d)$
 which follows along the lines of the inequality \cite[(1.3)]{van-Lint-Richert*65} by Van 
  Lint and Richert. Hence we have
  \begin{eqnarray*}
    \sum_{\substack{\e_1,\e_2}}
    \frac{|\lambda_{\e_1}(\q)\lambda_{\e_2}(\q)|}{\N[\e_1,\e_2]^{\alpha}}
    &\le&
      n^{18 n}z^{2(1-\alpha)}
      \sum_{\d, \atop{0< \N\d\le z}}\mu^2(\d)\frac{\phi_\alpha(\d)}{(2+\log(z/\N\d))^2\N\d^2} \\
      &\le&
     \frac{ n^{20 n}z^{2(1-\alpha)}}{(\log z)^2}
      \sum_{\d, \atop{0< \N\d\le z^{\alpha}}}
      \mu^2(\d)\frac{\phi_\alpha(\d)}{\N\d^2}
      ~+~ \frac{n^{18 n}z^{2(1-\alpha)}}{4} \sum_{\d, \atop{z^{\alpha}<\N\d\le z}} 
      \mu^2(\d)\frac{\phi_\alpha(\d)}{\N\d^2}
      \Bigl(\frac{\N\d}{z^\alpha}\Bigr)^{\frac{1-\alpha}{2}} \\
      &\le&
      n^{20 n}\frac{z^{2(1-\alpha)}}{(\log z)^2} \zeta_\K(2-\alpha)
      ~+~ n^{18 n}\frac{z^{\frac{(4-\alpha)}{2}(1-\alpha)}}{4}
      \zeta_\K\Bigl(1+\frac{1-\alpha}{2}\Bigr).
    \end{eqnarray*}
    We simplify the upper bound as follows;
    \begin{equation*}
      n^{20 n}\frac{z^{2/n}}{(\log z)^2}\biggl( (n+1)^{n}
      ~+~
      \frac{(\log z)^2}{4 z^{\frac{n-1}{2n^2}}}(2n+1)^{n}
      \biggr)
    \end{equation*}
    which is bounded above by
    \begin{equation*}
      n^{20n}\frac{z^{2/n}}{(\log z)^2}\biggl( (n+1)^{n}
      ~+~
      \frac{(\frac{4n^2}{n-1}\log z^{\frac{n-1}{4n^2}})^2}{4 z^{\frac{n-1}{2n^2}}}(2n+1)^{n}
    \biggr)
    ~\le~
    n^{27 n}\frac{z^{2/n}}{(\log
      z)^2}.
    \end{equation*}
%%%%%%%%%%%%%%%%%%%%%%%%%%%%%%%%%%%%%%%%%%%%%%%%
  Let us resume our study of $T_1$. The estimates above lead to
  \begin{multline*}
    \frac{h_{\K,\q}}{ X}T_1\le
    \frac{\alpha_\K \phi(\q) }{\N\q G_\q(z)}
    +
    36 \cdot 2^{14n} n^{8n}R_\K
    F(\q)
    \frac{z^2h_{\K,\q}}{\N\q }
    \frac{\N\q/X}{(2+\log z)^2}
    \\+~
    n \frac{h_{\K,\q}  }{\N\q }
    \frac{z\N\q}{X}
    ~+~
    n^{27 n} E(\K) F(\q)^{\frac{1}{n}}\log(3F(\q))^{n}
    \frac{h_{\K,\q}}{\N\q (\log z)^2}
    \left( \frac{z^2\N\q} {X}  \right)^{\frac{1}{n} },
  \end{multline*}
  where $E(\K)=1000
  n^{ 12n^2 }R_\K^{1/n}
  \bigl[\log\bigl((2n)^{4n}R_{\K} \bigr)\bigr]^{n}$ and $F(\q)=2^{r_1}\phi(\q)h_\K/h_{\K,\q}$ are defined in Theorem~\ref{asymfinal}. We employ Lemma~\ref{Gzbound} to bound $G_\q(z)$ from below, on assuming
  that $z\ge (10^6n)^{4n}|d_\K|^3$, getting
   \begin{multline*}
    \frac{h_{\K,\q}}{ X}T_1
    ~\le~
    \frac{1}{\log\frac{z}{e^4 n^4\sqrt{|d_{\K}|}}}
    \biggl(1
    +
    36 \cdot 2^{14n} n^{8n}R_\K
    2^{r_1}h_\K 
    \frac{z^2\N\q /X}{2+\log z}
    \\+~
    n 2^{r_1}h_\K
    \frac{z\N\q\log z}{X}
    ~+~
   n^{27 n}  E(\K) \frac{F(\q)^{\frac{1}{n}}\log(3F(\q))^{n}}{F(\q)}
    \frac{2^{r_1}h_\K}{ \log z}
    \left( \frac{z^2\N\q} {X}  \right)^{\frac{1}{n} }
    \biggr).
  \end{multline*}
  We notice that $z\log z\le \frac23 z^2/\log z$ when $z > 1$ and that
  \begin{align}
    \frac{y^{\frac{1}{n}}\log(3y)^{n}}{y}
    &=\biggl(\frac{\log(3y)}{y^{\frac{1}{n}-\frac{1}{n^2}}}\biggr)^{n}
      =3^{1-\frac1n}\biggl(\frac{\log(3y)}{(3y)^{\frac{1}{n}-\frac{1}{n^2}}}\biggr)^{n}
      \notag
    \\&=3^{1-\frac1n}\frac{n^{2n}}{(n-1)^n}
    \biggl(\frac{\log((3y)^{\frac1n-\frac1{n^2}})}{(3y)^{\frac{1}{n}-\frac{1}{n^2}}}\biggr)^{n}
    ~\le~ 3\frac{n^{2n}}{e^n(n-1)^n} ~\le~ n^{2n}.
    \label{simpl}
  \end{align}
  On assuming that $z^2\le X/\N\q$, we reach the upper bound
\begin{align*}
  \frac{h_{\K,\q}}{ X}T_1
  &\le~
    \frac{1}{\log\frac{z}{n^4\sqrt{|d_{\K}|}}}
    \biggl(1
    +
    \Bigl(36 n^{23 n}R_\K 
    +
    n^{2n}
    +
    E(\K)
    n^{30 n}
    \Bigr)\frac{h_\K }{ \log z}  \left( \frac{z^2\N\q} {X}  \right)^{\frac{1}{n} }\biggr)
  \\&\le~
  \frac{1}{\log\frac{z}{n^4\sqrt{|d_{\K}|}}}
  \biggl(1
    +
    \Bigl( y 
    +
    1
    +
    1000
  y^{1/n}
  (\log y)^{n}
    \Bigr)\frac{n^{ 27n^2 }h_\K}{\log z} \left( \frac{z^2\N\q} {X}  \right)^{\frac{1}{n} }\biggr)
  \end{align*}
  with $y=(2n)^{4n}R_\K$. Recall that $R_\K\ge 1/5$ by
  \cite{Friedman*89}, so that $y\ge13\,000$. A reasoning very similar to the
  one that led to inequality~\eqref{simpl}
  applies, getting
  \begin{equation*}
    \frac{h_{\K,\q}}{ X}T_1
    ~\le~
    \frac{1}{\log\frac{z}{n^4\sqrt{|d_{\K}|}}}
  \Bigl(1
    +
    \frac{n^{ 35n^2 }R_\K h_\K}{ \log z}  \left( \frac{z^2\N\q} {X}  \right)^{\frac{1}{n} }\Bigr).
      \end{equation*}
  Choosing $z^2 = \frac{X}{\N\q} (n^{ 35n^2 }R_\K h_\K)^{-n}$, we get
   \begin{equation*}
    \frac{h_{\K,\q}}{ X}T_1
    ~\le~
    \frac{1}{\log\frac{z}{n^4\sqrt{|d_{\K}|}}- 1}
    ~\le~
    \frac{2}{\log \frac{X}{\N\q} - \log \left(  n^8 |d_{\K}| n^{ 35n^3 } (R_\K h_\K)^n \right)}
  \end{equation*}
  and this completes the proof of our theorem.
\end{proof}
%%%%%%%%%% 

%%%%%%%%%%%%%%%%%%%%%%%%%%%%%%%%%%%%%%%%%%%%%%%%%%%%%%%%%%%%%%%%%%%%%%%%%% 
%%%%%%%%%%%%%%%%%%%%%%%%%%%%%%%%%%%%%%%%%%%%%%%%%%%%%%%%%%%%%%%%%%%%%%%%%%
\section{Finding enough primes}
%%%%%%%%%%%%%%%%%%%%%%%%%%%%%%%%%%%%%%%%%%%%%%%%%%%%%%%%%%%%%%%%%%%%%%%%%%
%%%%%%%%%%%%%%%%%%%%%%%%%%%%%%%%%%%%%%%%%%%%%%%%%%%%%%%%%%%%%%%%%%%%%%%%%%

Winckler in his PhD thesis (see Theorem 1.7 of \cite{Winckler*15})  proved an explicit version 
of Tchebotarev density Theorem. The case $L=\K$ provides us with the following explicit
version of Landau's prime ideal theorem.
%%%%%%%%%%%%%%%%%%%%%%%%%%
\begin{thm}[Winckler]
  \label{LandauPrimeIdealTheorem}
  For every $x \ge \exp (110000n_\K (\log(9|d_\K|^8 ))^2 )$, we have
  \begin{equation*}
    \sum_{\N\P\le x}1 
    ~=~ 
    \li(x) ~+~ \rO^*\bigl(\li(x^\beta)\bigr) ~+~
    \rO^*\biggl(10^{14}x\exp\frac{-\sqrt{\log x}}{12}\biggr),
  \end{equation*}
  where $\beta$ is the possible largest real zero of $\zeta_\K$ and $\li(x)=\int_2^x  \frac{dt}{\log t}$ is the
 usual  \emph{logarithmic integral function}.
\end{thm}
%%%%%%%%%%%%%%%%%%%%%%%%%%
This possible exceptional zero $\beta$ can be controlled by
the next lemma.
%%%%%%%%%%%
\begin{lem}[Kadiri and Wong \cite{KW}, 2022]
  \label{AKzero}
Any positive real zero $\rho$ of $\zeta_{\K}(s)$ satisfies
$1 - \rho \ge {1}/{|d_{\K}|^{12}}.$
\end{lem}
%%%%%%%%%%%

We need a lower bound for the number of prime ideals of degree
one, meaning we need to remove the primes of degree $>1$ from the
estimate of Theorem~\ref{LandauPrimeIdealTheorem}. This, and more, is
achieved by using the next Lemma. 
%%%%%%%%%%%%
\begin{lem}
  \label{RemoveParasites}
  We have
  \begin{equation*}
    \sum_{\substack{\P, k,\\ \N\P^k>p,\\ \N\P\le X}}1
    ~\le~
    3.64~n_\K\sqrt{X}.
  \end{equation*}
\end{lem}
%%%%%%%%%%%%
The summation is over powers $k \ge2$ for primes of any degree, or of
powers $k \ge1$ for primes of degree~$>1$. We denote by $p$ the rational
prime below~$\P$.
%%%%%%%%
\begin{proof}
  Each prime $p$ has at most $n_\K$ prime ideals above itself. For
  each such prime, say of norm $p^f$, the contribution to the sum is
  at most $\frac{\log X}{f\log p}\le (\log X)/\log 2$ and a prime $p$ that contributes verifies
  $p\le \sqrt{X}$. We use the bound for $\pi(x)\le 1.26x/\log x$ from
  \cite[Corollary 1]{Rosser-Schoenfeld*62} by Rosser and Schoenfeld valid for
  $x=\sqrt{X}>1$. The lemma follows. 
\end{proof}
%%%%%%%%

%%%%%%%%%
\begin{lem}
  \label{lowerbound}
  When $x\ge \exp(|d_\K|^{30})$, we have
 $\displaystyle
    \sumflat_{\substack{\P \atop{\N\P\le x}}}1
    ~\ge~
    \frac{x}{\log x} ~+~ \frac{(1-10^{-10})x}{(\log x)^2}$.
\end{lem}
%%%%%%%%%

%%%%%%%%%
\begin{proof}
  Let $\pi^\flat_\K(x)$ be the quantity to estimate.
  We start by combining Theorem~\ref{LandauPrimeIdealTheorem} together
  with
  Lemma~\ref{RemoveParasites} to obtain
  \begin{equation*}
    \pi^\flat_\K(x)
    \ge
    \li(x)
    -\li(x^\beta)
    -10^{14}x\exp\biggl(\frac{-\sqrt{\log x}}{12}\biggr)
    -3.64n_\K\sqrt{x}.
  \end{equation*}
  By Theorem~\ref{AKzero}, we have $\beta\le 1-|d_\K|^{-115}$, while,
  by Lemma~\ref{rootdisc}, we have $n_\K\le 3\log|d_\K|$. This gives
  us a lower bound solely in terms of the discriminant. We also need
  to handle the logarithmic integral. Using integration by parts, we find
  that
  \begin{equation*}
    \li(y)=\frac{y}{\log y}-\frac{2}{\log 2}
    +\int_2^{y}\frac{dt}{\log^2t}
    ~\ge~ \frac{y}{\log y}-\frac{2}{\log
      2}+\frac{\li(y)}{\log y}.
  \end{equation*}
  We first deduce from this that $\li(y)\ge y/\log y$ when
  $y\ge 7.5$, by neglecting the additional term $\li(y)/\log y$.
  On incorporating it, we get $\li(y)\ge \frac{y}{\log
    y}+\frac{y}{(\log y)^2}$.
  As for an upper bound, we work rather trivially:
  \begin{equation*}
    \li(y)
   ~ \le~ \int_2^{\sqrt{y}}\frac{dt}{\log 2}+\int_{\sqrt{y}}^y\frac{dt}{\log\sqrt{y}}
    ~\le~ \frac{2y}{\log y}+\frac{\sqrt{y}}{\log 2}
    ~\le~ \frac{3\,y}{\log y}
  \end{equation*}
  when $y\ge 20$. This leads to, with $d=|d_\K|$,
  \begin{align*}
    \biggl(\frac{\log x}{x}\pi^\flat_\K(x)-1\biggr)\log x
    &\ge~
      1-\frac{3e^{-\frac{\log x}{d^{12}}}\log x}{\beta}
      -10^{14}(\log x)^2e^{\frac{-\sqrt{\log x}}{12}}
      -12(\log d)\frac{(\log x)^2}{\sqrt{x}},
    \\&\ge~
    1-6d^{30}e^{-d^{18}}
    -10^{14} d^{60} e^{\frac{-d^{15}}{12}}
    -12(\log d)d^{60}e^{-\frac{d^{30}}{2}}
    \\&\ge~ 1-10^{-10}
  \end{align*}
  since $d\ge2$. This completes the proof of the lemma.
\end{proof}
%%%%%%%%%

%%%%%%%%% 
\begin{lem}
  \label{wlowerbound}
  When $x\ge \exp(|d_\K|^{30})$, we have
  $\displaystyle
    \sumflat_{\substack{\P}}w_0(\N\P/x)
    \ge
    \frac{x\|w_0\|_1}{\log x}
    +\frac{x\|w_0\|_1}{5(\log x)^2}$.
\end{lem}
%%%%%%%%%

%%%%%%%%%
\begin{proof}
  Since $w_0$ is bounded above by~1 and has support within $[0,1]$, we
  may again use Lemma~\ref{RemoveParasites} to handle the condition
  `$\P$ of degree~1'. Whence, on denoting $T(w_0)$ the sum to be studied, we find that
  \begin{equation*}
    T(w_0)
    \ge \sum_{\substack{\P}}w_0(\N\P/x)-3.64n_\K\sqrt{x}.
  \end{equation*}
  As $w_0(y)= - \int_y^1w_0'(t)dt$, we get
  \begin{equation*}
    T(w_0)\ge
   - \int_0^1
    \bigl(\sum_{\N\P\le tx}1\bigr) w_0'(t)dt-3.64n_\K\sqrt{x}.
  \end{equation*}
  Of course $w_0'(t)=0$ when $t\le 1/10$, so we may assume that $tx\ge x/10$
  which is thus larger than $ \exp (110000n_\K (\log(9|d_\K|^8 ))^2
  )$. Theorem~\ref{LandauPrimeIdealTheorem} applies and yields
  \begin{align*}
    T(w_0)
    &\ge~
   - \int_0^1
    \li(xt) w_0'(t)dt
    -\|w_0'\|_1\biggl(\li(x^\beta)+10^{14}x\exp\frac{-\sqrt{\log x}}{12}\biggr)
    -3.64n_\K\sqrt{x}
    \\&\ge~
    \int_0^1
     \frac{xw_0(t)dt}{\log(xt)}
    -\|w_0'\|_1\biggl(\li(x^\beta)+10^{14}x\exp\frac{-\sqrt{\log x}}{12}\biggr)
    -3.64n_\K\sqrt{x}
    \\&\ge~
    \frac{x\|w_0\|_1}{\log x}
    +\int_0^1\frac{xw_0(t)\log(1/t)dt}{(\log x)\log(x t)}
    -\|w_0'\|_1\biggl(\li(x^\beta)+10^{14}x\exp\frac{-\sqrt{\log x}}{12}\biggr)
    -3.64n_\K\sqrt{x}.
  \end{align*}
  All of that is valid for a rather general non-negative
  function~$w$. When it comes to $w=w_0$, we have
  $10\sqrt{n_{\K}}\|w_0\|_1\in[2, 15]$ by
  Lemma~\ref{studyw0} and $\|w_0'\|_1=2$. Further,
by applying Lemma~\ref{rootdisc}, we get $\sqrt{n_{\K}} \le
  \sqrt{(\log|d_\K|)/\log(\pi/2)}$
  which is not more than $|d_\K|^2$. We finally notice that
  \begin{equation*}
    \|w_0\log(1/t)\|_1
    ~\ge~ \int_{1/10}^{11/20}w_0(t)\log(1/t)dt
    ~\ge~ \log(20/11)\frac{\|w_0\|_1}2 
    ~\ge~ \frac{\|w_0\|_1}4.
  \end{equation*}
  Hence, and following estimates very similar to the ones done during
  the proof of Lemma~\ref{lowerbound}, we find that
  \begin{align*}
    T(w_0)
    &\ge~
    \frac{x\|w_0\|_1}{\log x}
    +\frac{x\|w_0\|_1}{4(\log x)^2}\biggl(
    1
    - 832 \sqrt{\log |d_\K|} x^{\beta-1}\log x
    - 757  (\log |d_\K|)^{3/2}x^{-1/2}(\log x)^2
    -10^{-6} \frac{\sqrt{\log |d_\K|}}{\log x} 
    \biggr)
    \\&\ge~
    \frac{x\|w_0\|_1}{\log x}
    +\frac{x\|w_0\|_1}{5(\log x)^2}
  \end{align*}
  as required.
\end{proof}
%%%%%%%%%

%%%%%%%%%%%%%%%%%%%%%%%%%%%%%%%%%%%%%%%%%%%%%%%%%%
%%%%%%%%%%%%%%%%%%%%%%%%%%%%%%%%%%%%%%%%%%%%%%%%%%
\section{Proof of \thmref{mainthm}}
%%%%%%%%%%%%%%%%%%%%%%%%%%%%%%%%%%%%%%%%%%%%%%%%%%
%%%%%%%%%%%%%%%%%%%%%%%%%%%%%%%%%%%%%%%%%%%%%%%%%%
In order to prove \thmref{mainthm}, we need the following lemmas.
%%%%%%%%%%
\begin{lem}
  \label{simplifytK}
 We have $n_\K^{ 48n_\K^3 }(R_\K h_\K)^{n_{\K}} \ge 
 10^{25n_\K} n_{\K}^{7n_{\K}}$. 
\end{lem}
%%%%%%%%%%

%%%%%%%%%%
\begin{proof}
By ~\cite{Friedman*89}, we have $R_\K/ |\mu_\K| \ge 9/100$ which implies $R_\K\ge 9/100$. 
It is thus enough to check the inequality
  \begin{equation*}
    n_\K^{ 48n_\K^3 } \left( \frac{9}{100} \right)^{n_{\K}} ~\ge~ 
    10^{25n_\K} n_{\K}^{7n_{\K}} 
    \end{equation*}
  which is readily seen to hold true as $n_{\K } \ge 2$.
\end{proof}
%%%%%%%%%%

% %%%%%%%%%%

\begin{lem}[Kneser's Theorem \cite{Kneser*53}, 1953]\label{KT}
Let $G$ be a finite abelian group and $\B$ be a non-empty
subset of $G$. Also let 
$
H = \{ g \in G  ~|~ g + \B + \B = \B + \B \}
$ 
be the stabiliser of the set $\B + \B$. If $\B$ intersects $\lambda$ many cosets 
of $H$, then
$$
|\B+\B| \geq (2\lambda -1 )|H|.
$$
\end{lem}

%%%%%%%%%%%%%%%
\begin{proof}[Proof of Theorem~\ref{mainthm}]
For any un-ramified prime ideal $\P$ of degree one, we denote $\N(\P)$ by $p$.
Also let $\A$ be the subset of $H_{\q}(\K)$ defined~by
\begin{equation}
\A = \{ [\a] \in H_{\q}(\K) ~|~ \exists ~\P \text{ with } \N(\P)=p < X, [\P]=[\a]\}.
\end{equation}
Let us set $u(\K)= n_{\K}^{ 48n_{\K}^3 } |d_\K|^6(R_\K h_\K)^{n_{\K}}$. When
$X>u(\K)\N\q$, Theorem~\ref{bt-tri} gives us that
\begin{equation*}
  \sumflat_{\substack{\N\P\le X}}1
  ~\le~
  \frac{2 |\A| X}
  {h_{\K,\q}\log\frac{X}{u(\K)\N\q}}.
\end{equation*}
On the other side and when $X\ge \exp(|d_\K|^{30})$,
Lemma~\ref{lowerbound} ensures us that
\begin{equation*}
  \frac{X}{2(\log X)^2}+\frac{X}{\log X}~\le~ \sumflat_{\substack{\N\P\le X}}1.
\end{equation*}
On assuming that the required conditions hold true, a comparison of
both inequalities gives us
\begin{equation*}
   \frac{|\A|}
   {h_{\K,\q}}
   ~\ge~
   \frac{\log\frac{X}{u(\K)\N\q}}{2\log X}
   ~+~
   \frac{\log\frac{X}{u(\K)\N\q}}{4(\log X)^2}.
\end{equation*}
Take $X=(t(\K)\N\q)^3$, where $t(\K)$ is defined in \eqref{deftK}. The above inequality gives us
\begin{equation}
  \label{inieq}
  \frac{|\A|} {h_{\K,\q}}
  ~\ge~
  \frac13 
  ~+~
  \frac{1}{18\log
    t(\K)+ 18\log\N\q}.
\end{equation}
From now onwards, let us set $G =H_\q(\K)$.
Also let $H$ be the stabilizer of $\A\cdot \A$ in $G$, $y$ be the index of
$H$ in $G$ and $\lambda$ be the number of cosets of $H$ that intersect
$\A$. 

By Kneser's Theorem, we have
\begin{equation*}
  |\A \cdot \A|
  ~\ge~
   (2\lambda-1)|H| 
   ~=~ \frac{2\lambda-1}{y}|G|.
\end{equation*}
We further know that
\begin{equation*}
 \lambda ~\ge~ \Bigl\lceil{\frac{|\A|}{|H|}\Bigr\rceil} 
\end{equation*}
and that $\lambda$ is an integer.
Furthermore, to be sure that $\A \cdot \A \cdot\A =G$, we only
need
\begin{equation}\label{eq:8}
 |\A|  + \frac{2\lambda-1}{y}|G|  ~>~  |G|
 \phantom{m}\text{i.e.}\phantom{m}
 \frac{|\A|}{|G|}+\frac{2\lambda-1}{y}  ~>~ 1.
\end{equation}
This follows from the following observation.
Given any $[\b] \in G$, we consider the
set 
$$
[\b]\A^{-1} = \{[\b] [\a] ~~:~ [\a^{-1}] \in \A\}.
$$
For any $[\b] \in G$ if  $[\b]\A^{-1} \cap \A\cdot\A \neq \emptyset$,
then $[\b] \in \A\cdot\A\cdot\A$. However $|[\b]\A^{-1}| = |\A|$.
Therefore by Pigeon-hole principle if $|\A| + |\A \cdot \A| > |G|$, then
$\A \cdot \A\cdot \A=G$.

From the lower bound $\lambda \ge |\A|/|H|$, we observe
that it suffices to show that
\begin{equation} \label{eq:9}
 3 \frac{|\A|}{|G|} - \frac{1}{y} ~>~ 1.
\end{equation}
Let us discuss the possible values of $y$.
%%%%%%%%%%%%%%%%%%%%%%%%%%%%
\subsubsection*{$\blacksquare$ Large values of $y$}
%%%%%%%%%%%%%%%%%%%%%%%%%%%%
By~\eqref{inieq}, the above inequality \eqref{eq:9} will be satisfied if we have
$$
\frac{1}{9\log t(\K)~+~ 9\log \N\q}-\frac1y ~>~ 0
$$
This implies that the inequality \eqref{eq:9} holds when
$y > 9 \log t(\K) + 9\log \N\q$.
%%%%%%%%%%%%%%%%%%%%%%%%%%%%%%%%%%%%%%%%%%%%%%
\subsubsection*{$\blacksquare$ The case $y=1$}
%%%%%%%%%%%%%%%%%%%%%%%%%%%%%%%%%%%%%%%%%%%%%%
This is the case when $H=G$ and
thus $\A \cdot \A=G$. Hence $\A \cdot \A \cdot\A =G$.
%%%%%%%%%%%%%%%%%%%%%%%%%%%%%%%%%%%%%%%%%%%%%%%
\subsubsection*{$\blacksquare$ The case $y=2$}
%%%%%%%%%%%%%%%%%%%%%%%%%%%%%%%%%%%%%%%%%%%%%%%
So $H$ is a quadratic subgroup of $G$. Using Lemma~\ref{simplifytK}, 
we have
$$
X
~\ge~ ( 10^{25n_\K} n_{\K}^{7n_{\K}} |d_\K|^{4/3} \N\q)^3 
~\ge~ 10^{25n_\K} n_{\K}^{7n_{\K}} |d_\K|^{4}  \N\q^3.
$$
Theorem~\ref{degreeoneprime} tells us that the subgroup
generated by $\A$ is $G$. It follows that 
$\A$ contains an element of $G\setminus H$.  By \thmref{primeinkernel},
we have that the subset $\A$ also contains an element of $H$. Indeed, we
have $X\ge 8(10^{31} n_\K^7)^{n_\K}  |d_\K|^4 \N\q^2$.
 As $\A\cdot\A$ is a union of cosets mod~$H$, we conclude
 that $ \A\cdot\A=G$, whence $\A\cdot\A\cdot\A=G$.

%%%%%%%%%%%%%%%%%%%%%%%%%%%%%%%%%%%%%%%%%%%%%%%
\subsubsection*{$\blacksquare$ Medium values of $y\nequiv2 \bmod{3}$}
%%%%%%%%%%%%%%%%%%%%%%%%%%%%%%%%%%%%%%%%%%%%%%%

Since $\A\cdot \A$ is a union of cosets modulo $H$, it is enough to
check that $\A/H\cdot(\A\cdot \A/H)$ covers $G/H$. This means that it
is enough to assume that $\A$ is a union of ($\lambda$ many) cosets modulo~$H$, from
which we infer that we only need to prove that
$\frac{\lambda}{y}+\frac{2\lambda-1}{y}> 1$, i.e. 
$\lambda>(y+1)/3$. Note that \eqref{inieq} implies that
$\lambda>y/3$. When $y\equiv 0 \bmod{3}$, then $\lambda \ge y/3 +1$.
When $y\equiv 1 \bmod{3}$, the integer part of $y/3$ is at
least $(y+2)/3$, so that the inequality $\lambda>(y+1)/3$ is
satisfied. There remains the values of $y\ge3$ that are $2 \bmod{3}$
and below $9 \log t(\K)+ 9\log \N\q$.

%%%%%%%%%%%%%%%%%%%%%%%%%%%%%%%%%%%%%%%%%%%%%%%
\subsubsection*{$\blacksquare$ Medium values of $y\equiv2 \bmod{3}$ and $y\ge5$}
%%%%%%%%%%%%%%%%%%%%%%%%%%%%%%%%%%%%%%%%%%%%%%%
We use \thmref{bt} for $y \le 9 \log t(\K)+ 9 \log \N\q$ together with
Lemma~\ref{wlowerbound}. 
This gives us, with $n_\K=n$,
\begin{align*}
  \frac{X\|w_0\|_1}{\log X}
  +\frac{X\|w_0\|_1}{5(\log X)^2}
  &\le~
  \sumflat_{\substack{\P}}w_0(\N\P/X)
  ~\le~
  \sum_{ [\a]H \in\A/H}\sum_{\P \in [\a]H}  w_0\left(\N\P/X\right)
  \\&\le~
  \lambda
  \frac{2 \|w_0\|_1 X}
  {y\log\frac{X\|w_0\|_1/[20000 (\|w_0^{(n+3)}\|_\infty
  +5\|w_0\|_1)\log(|d_\K|\N\q)^{n}]}
    {2^{22n}y\sqrt{|d_{\K}|^3\N(\q)}}}.
\end{align*}
Lemma~\ref{studyw0} gives us
\begin{equation*}
  \frac{\|w_0\|_1}{\|w_0^{(n+3)}\|_\infty +5\|w_0\|_1}
  ~\ge~ \frac{2\sqrt{n}}{ (40n)^{n+4} + 75 \sqrt{n} }
  ~\ge~ \frac{\sqrt{n}}{(40n)^{n +4}}.
\end{equation*}
We have thus obtained the inequality
\begin{equation*}
  \log\frac{\sqrt{n}}{20000(40n)^{n+4}
    2^{22n}|d_{\K}|^{3/2}}
  \frac{X}
  {9\log (t(\K)\N\q)\sqrt{\N(\q)}\log(|d_\K|\N\q)^{n}}
  ~\le~ \frac{2\lambda}{y}\log X
\end{equation*}
i.e. also
\begin{equation*}
  1-\frac{1}{\log X}
  \log\biggl(180000\cdot 2^{22n}(40n)^{n+4}|d_{\K}|^{3/2}
  {\log (t(\K)\N\q)\sqrt{\N(\q)}\log(|d_\K|\N\q)^{n}}
  \biggr)
  ~\le~ \frac{2\lambda}{y}.
\end{equation*}
We readily check that $180000\cdot 2^{22n}(40n)^{n+4} |d_{\K}|^{3/2} \le \sqrt{t(\K)}$, so that the
    above inequality implies that
    \begin{equation*}
      1-\frac1{3\log V}\log(\sqrt{V}(\log V)^{1+n}) ~\le~ \frac{2\lambda}{y},
    \end{equation*}
    where $V=t(\K)\N\q$,
    while we are to prove that $\lambda > (y+1)/3$ (see the discussion of
    the case ``Medium values of $y\nequiv2 \bmod{3}$"). We should thus prove
    that $\log(\sqrt{V}(\log V)^{1+n})<\frac{3}{5}\log V$, i.e.
    \begin{equation*}
      (\log V)^{10 + 10 n} ~<~ V.
    \end{equation*}
    Again by using Lemma~\ref{rootdisc}, we see that
    $n\le \frac{1}{30\log(\pi/2)}\log\log
    t(\K)\le \frac{1}{13}\log\log V$.
    With $W=\log V$, we should thus prove that
    \begin{equation*}
      10\log W ~+~ \frac{10}{13} (\log W)^2 ~<~ W.
    \end{equation*}
    This inequality happens to be always satisfied, which
    concludes our proof.
\end{proof}
%%%%%%%%%%%%

\noindent
{\bf Acknowledgements.} 
Research of this article was partially supported by Indo-French Program in Mathematics (IFPM).  
All authors would like to thank IFPM for financial support. The first, second and fourth authors 
acknowledge the support of SPARC project~445.
The second author would also like to thank MTR/2018/000201 and DAE number 
theory plan project for partial financial support.

\end{document}